\documentclass{dcds-b}
\usepackage{amsmath}
  \usepackage{paralist}
  \usepackage{graphics} 
  \usepackage{epsfig} 
 \usepackage[colorlinks=true]{hyperref}
 \hypersetup{urlcolor=blue, citecolor=red}

\def\currentvolume{XX}
 \def\currentissue{X}
  \def\currentyear{2010}
   \def\currentmonth{Month X}
    \def\ppages{1--x}
     \def\DOI{yy.yyyy/dcdsb.2010.yy.y}

\DeclareMathOperator*{\esssup}{ess\,sup}
\newcommand{\ra}{\rightarrow}
\newcommand{\RR}{\mathbb{R}}
\newcommand{\Sn}{\mathbb{S}^{n-1}}
\newcommand{\Se}{\mathbb{S}^1}

\newcommand{\RRn}{\RR^n}

\newcommand{\BV}{\mathcal{B}(V)}
\newcommand{\BS}{\mathcal{B}(\Sn)}

\newcommand{\X}{\mathbb{X}}
\newcommand{\A}{\mathcal{A}}
\newcommand{\G}{\mathcal{G}}
\newcommand{\D}{\mathcal{D}}
\newcommand{\U}{\mathcal{U}}

\newcommand{\LL}{\mathcal{L}}
\newcommand{\ud}{\mathrm{d}}

\newcommand{\e}{\varepsilon}

\newtheorem{theorem}{Theorem}[section]
\newtheorem{corollary}{Corollary}
\newtheorem{lemma}[theorem]{Lemma}
\theoremstyle{definition}
\newtheorem{definition}[theorem]{Definition}
\newtheorem{remark}{Remark}
\newtheorem{example}[theorem]{Example}

\usepackage{color,ulem}
\title[Cell Movement in Network Tissues]{Mathematical Analysis of a Kinetic
Model for Cell Movement in Network Tissues}
\author[T. Hillen, P. Hinow and Z. Wang]{}


\email{thillen@math.ualberta.ca}
\email{hinow@uwm.edu}
\email{zhian.wang@vanderbilt.edu}

\thanks{TH was partially supported by NSERC and MITACS}

\begin{document}
\maketitle

\centerline{\scshape Thomas Hillen }
{\footnotesize
 \centerline{Department of Mathematical and Statistical Sciences}
 \centerline{Centre for Mathematical Biology}
\centerline{University of Alberta, Edmonton, T6G 2G1, Canada}  
} 
\medskip
\centerline{\scshape Peter Hinow }
{\footnotesize
 \centerline{Department of Mathematical Sciences}
 \centerline{University of Wisconsin -- Milwaukee}
 \centerline{P.O. Box 413, Milwaukee, WI 53201-0413, USA}
} 
\medskip
\centerline{\scshape Zhi-An Wang}
{\footnotesize
 \centerline{ Department of Mathematics }
   \centerline{Vanderbilt University, Nashville, TN 37240, USA } 
}
\bigskip

 \centerline{(Communicated by Kevin Painter)}

\begin{abstract} 
Mesenchymal motion describes the movement of cells in biological tissues formed
by fibre networks. An important example is the migration of tumour cells through
collagen networks during the process of metastasis formation. We investigate the
mesenchymal motion model proposed by T.~Hillen in \cite{Hillen06} in higher
dimensions. We formulate the problem as an evolution equation in a Banach space
of measure-valued functions and use methods from semigroup theory to show the
global existence of classical solutions. We investigate steady states of the
model and show that patterns of network type exist as steady states. For the
case of constant fibre distribution, we find an explicit solution and we prove
the convergence to the parabolic limit.
\end{abstract}

\section{Introduction}\label{section:introduction}

\medskip
Friedl and collaborators \cite{Friedl} observed mesenchymal tumour  cells 
as they move in a field of collagen fibres and change their velocities according to the local orientation of the fibres. At the same time, the cells also remodel the fibres, primarily
through expression of matrix-degrading enzymes (proteases) that cut selected
fibres. In \cite{Hillen06}, the author introduced a mathematical model for this process of mesenchymal cell movement in fibrous tissues. Recent analysis of this and similar models 
 \cite{Hillen06, Painter, Chauviere1, Chauviere2, WHL} revealed the existence of biologically meaningful measure valued solutions, which correspond to tissue and cell alignment. Hence a sophisticated existence theory is needed. In this paper we will formulate the mesenchymal transport
model proposed in \cite{Hillen06} as a semilinear evolution equation in a Banach
space of measure-valued functions. We apply classical theory of semigroups of
operators and a Banach Fixed Point argument to show well-posedness of the
problem (Section \ref{global}). With the correct theoretical framework in place, we are then able to
classify possible steady states, whereby we introfduce a new notation of {\it pointwise steady states}, which are meant to resemble the network patterns which were observed numerically in \cite{Painter}. Moreover,  we rigorously study the parabolic limit (diffusion limit) of the kinetic model in the measure-valued context. We
show convergence to the diffusion limit for constant fibre distribution  in
Section \ref{limit}.

The existence theory here employs a mild solution formulation which is
based on a variation of constant formula. The solutions are functions in $L^1$
in space and measures in velocity. It turns out that this definition is too 
``weak'' in the sense that it does not provide a nice representation of the
global network patterns observed numerically. Hence here we introduce a
sub-class of steady states, which we call \textit{pointwise steady states}.
First of all, we show that pointwise steady states do exist. Secondly, pointwise
steady states allow for a representation of network patterns. Our results
include a classification of possible network intersections. 

In the model proposed in \cite{Hillen06}, undirected and directed tissues
were distinguished. In undirected tissues (e.g.~collagen), fibres are
symmetric and both directions are identical, a situation that somewhat
resembles a nematic liquid crystal \cite{Virga}. In directed tissues, fibres
are asymmetric and the two ends can be distinguished. From the mathematical
point of view, which we adopt in the present paper, both cases are completely
analogous. Hence we focus on the case of undirected tissues. We refer to 
\cite{Hillen06} for the  biological assumptions and the detailed mathematical
derivation of the model. 

The model studied here is specifically designed for mesenchymal cell movement
in network tissues via contact guidance and degradation of the extracellular
matrix (ECM). Painter \cite{Painter} has extended this model in various
directions. 
His model variations allow (i) to choose between amoeboid and mesenchymal motion, (ii) to 
place different weights between diffusive movement and movement by contact
guidance, (iii) to include ECM degradation as well as production, (iv) to
include ECM remodelling or lack thereof, (v) to study focussed
protease release at the cell tip versus unfocussed ECM degradation via a
diffusible proteolytic enzyme. 
Many of these modifications lead to the same pattern formation properties as
observed for the initial model. All of these modifications show the same
mathematical challenges, namely the description of aligned tissue as weak
solutions and orientation driven instabilities. Hence we believe that the
results which we present here are representative for a large class of kinetic
models for cell movement in tissues and they can be generalized to many other
cases.

In \cite{Hillen06}, the techniques of moment closure, parabolic and
hydrodynamic scaling were used to study the macroscopic limits of the system
that we later restate in equation \eqref{old_model}. The resulting
macroscopic models have the form of drift-diffusion equations where the mean
drift velocity  is given by the mean orientation of the tissue and the
diffusion  tensor is given  by the variance-covariance matrix of the tissue
orientations. Model  \eqref{old_model} has been extended in
\cite{Chauviere1,Chauviere2} to  include  cell-cell interactions and
chemotaxis. The corresponding diffusion limit was  formally obtained in
these papers.

In the case of chemotaxis, a system of a transport equation for the cell motion
coupled to a parabolic or elliptic equation for the chemical signal was studied
by Alt \cite{Alt80}, Chalub \textit{et al.} \cite{Chalub} and Hwang \textit{et
al.}
\cite{HKS1,HKS2}. Local and global existence of solutions were studied and the
macroscopic limits were proved rigorously in \cite{Chalub,HKS1,HKS2}. However,
these authors assumed the existence of an equilibrium velocity distribution for
cells that is in $L^{\infty}(V)$ where $V$ denotes the space of velocities. For
the mesenchymal motion model, it is necessary to allow for complete alignments
of either fibres or cells, corresponding to Dirac measures on $V$ or the space
of directions, the unit sphere $\Sn$. In particular, assumption (A0) in paper
\cite{Chalub} does not apply here and hence their respective results can not be
applied directly to the mesenchymal motion model. 

In Section \ref{s:statement} we formulate the model and we introduce suitable
function spaces and operators. Our first main result on global existence of
measure-valued solutions is given in Section \ref{s:existence}. In Section
\ref{steady_states} we present a definition and classification of pointwise
steady states. In Section \ref{characteristics} we assume that the fibre density
$q$ is a given function of $x, t$. In that case we find an explicit solution of
the kinetic equation using the methods of characteristics. If moreover, the
fibre distribution is constant in time and space, then we prove the
convergence to a parabolic limit.   It appears to be impossible to prove
convergence to the parabolic limit for arbitrary time- and space dependent fibre
distributions. This confirms numerical observations of Painter \cite{Painter},
who investigated the mesenchymal motion model and found interesting cases of
pattern formation of network type (see Figure \ref{fig_Kevin}). In the diffusion
limit, however, the patterns disappear in the numerical simulation. This
indicates that there is a significant difference in the asymptotics of the
kinetic model and the diffusion limit for timely varying tissue networks.  

\section{Formulation of the Problem}\label{s:statement}

\subsection{The Model}
\begin{figure}[ht]
\begin{center}
\includegraphics[width=8cm]{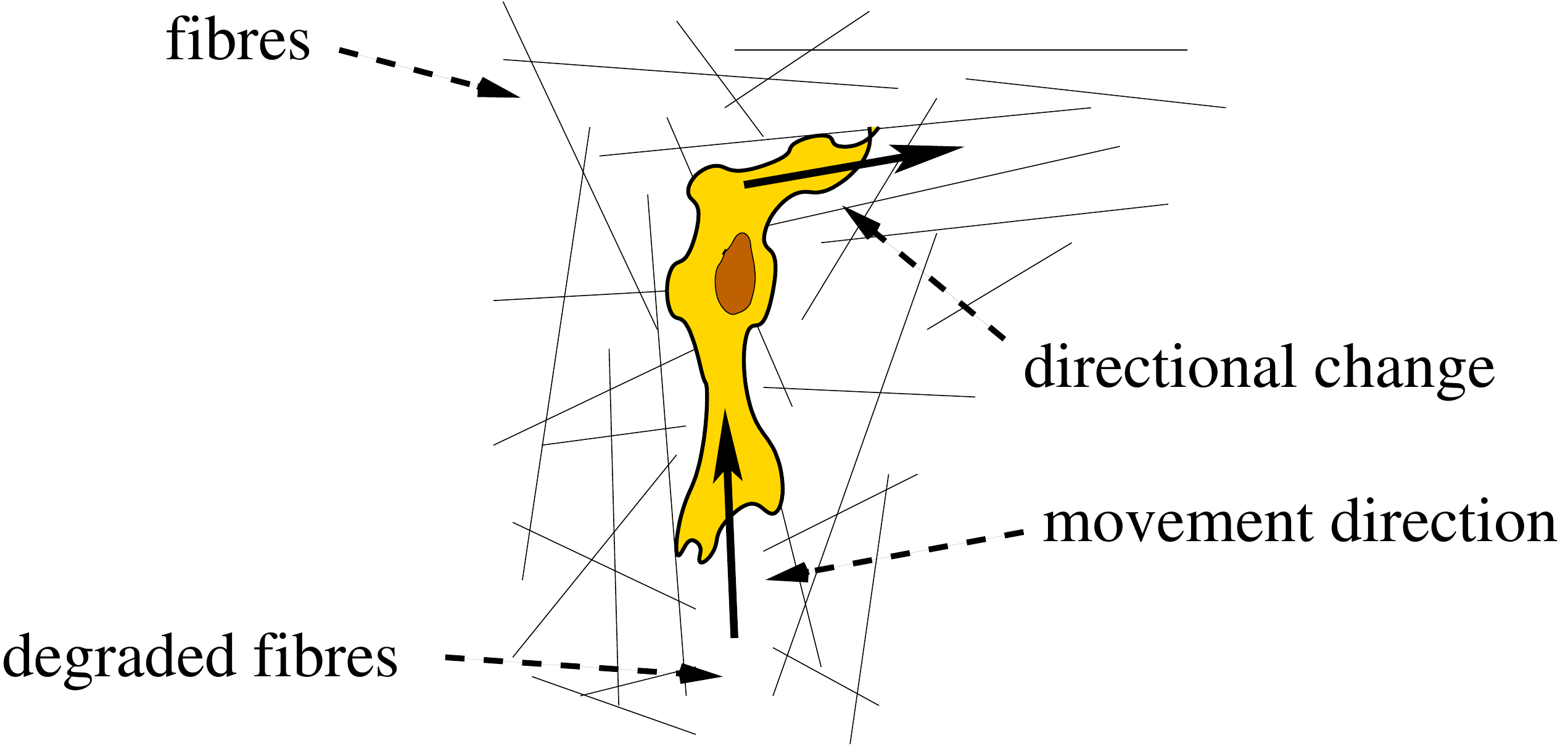}
\caption{Schematic of the model (\ref{old_model}) for cell movement in network
tissues, including directional changes, contact guidance and fibre
degradation.}\label{f:M5}
\end{center}
\end{figure}

We briefly recall the kinetic model for  mesenchymal motion from \cite{Hillen06}
for the undirected case. The distribution $p(x,t,v)$ describes the cell density
at time $t\ge0$, location $x\in\RRn$ and velocity $v\in V$. Throughout the paper
we assume that $V$ is a product $V=[s_1,s_2]\times\Sn$, where $0\le s_1\le
s_2<\infty$ is the range of possible speeds. If $s_1=s_2$ then we assume
$s_1>0$. The fibre network is described by the distribution $q(x,t,\theta)$ with
$\theta\in\Sn$, the $(n-1)$-dimensional unit sphere in $\RRn$. A schematic of
the model is given in Figure \ref{f:M5}. 

 The model for
mesenchymal motion from \cite{Hillen06} reads  
\begin{equation}\label{old_model}
\begin{aligned}
&\frac{\partial p(x,t,v)}{\partial t} +v\cdot \nabla p(x,t,v)  =-\mu
p(x,t,v)+\mu\bar{p}(x,t)\,\tilde{q}(x,t,v), \\
&\frac{\partial q(x,t,\theta)}{\partial t} =
\kappa(\Pi_u(p(x,t,v))-A_u(p(x,t,v),q(x,t,\theta)))\bar{p}(x,t) q(x,t,\theta),
\\
&p(x,0)=p_0(x), \quad q(x,0)=q_0(x),
\end{aligned}
\end{equation}
where $\mu$ and $\kappa$ are positive constants.  The transport term $v\cdot
\nabla p$ indicates that cells move with their velocity.  
 The right hand side of the first equation describes the reorientation of the cells in the field of
fibres. Turning away from their old direction at rate $\mu$, they turn into a
new direction with a probability that corresponds to the fibre distribution $q$.
The new speed is chosen from the interval $[s_1,s_2]$. The cells degrade (at
rate $\kappa$) those fibres that they meet at an approximately right angle while
they leave fibres that are parallel to their own orientation unchanged. The
exact definitions of the corresponding terms in system
\eqref{old_model} requires some mathematical details that are given in the
next section. The expressions $\bar{p}$, $\tilde{q}$, $\Pi_u(p)$ and $A_u(p,q)$
are defined in equations \eqref{mass}, \eqref{lifting}, \eqref{mean_projection}
and \eqref{relative_alignment}, respectively.

Painter showed in \cite{Painter} that the second equation of
\eqref{old_model} arises if instead of ECM degradation one assumes that the
cells realign the tissue. This would be the case for fibroblasts, who do remodel
the fibre newtork without destroying it. In that case the term $-A_u$ measures
the fibre degradation while $\Pi_u$ describes the fibre production such that the
total amount of fibre mass is preserved.

\subsection{Spaces and Operators}\label{spaces}
We show  in Section \ref{steady_states} that Dirac measures occur as  meaningful
steady states. Hence we need to construct a solution framework that allows for
measure-valued solutions. Let $\Omega=\RRn$ be the spatial domain in which
particles are able to move. 

Let $\BV$ denote the space of regular signed real-valued (finite) Borel measures
on $V$. For $p\in\BV$ let $p=p^+-p^-$ be its Hahn-Jordan decomposition and 
$|p|=p^++p^-$ its variation \cite{Cohn}. When equipped with the total variation
norm (the following notations are used interchangeably throughout the paper)
\begin{equation*}
||p||_{\BV} = |p|(V) = \int_V \ud |p|(v) = \int_V |p|(\ud v),
\end{equation*}
$\BV$ is a Banach space  \cite[Proposition 4.1.7]{Cohn}. Analogously, $\BS$ will
denote the Banach space of regular signed Borel measures on $\Sn$ equipped with
the total variation norm. Naturally, we are interested in solutions
taking values among non-negative measures only.  Let 
\begin{equation*}
\begin{aligned}
\X_1&=L^1(\RRn,\BV), \\ 
\X_2&=L^\infty(\RRn,\BS), \\
\X&=\X_1 \times \X_2,
\end{aligned}
\end{equation*}
equipped with norms \begin{equation*}
\begin{aligned}
||p||_{\X_1}&=\int_{\RRn} ||p(x)||_{\BV}\,\ud x,\\ ||q||_{\X_2}&=\esssup_{x\in \RRn} ||q(x)||_{\BS}, \\
||(p,q)||_{\X}&=||p||_{\X_1}+||q||_{\X_2}.
\end{aligned}
\end{equation*}
We denote the positive cones of the spaces $\X_1,\,\X_2$ and $\X$ by
$\X_1^+,\,\X_2^+$ and $\X^+$, respectively. We will write
\begin{equation*}
||p||_\infty = \esssup_{x\in \RRn} ||p(x)||_{\BV}
\end{equation*}
for those $p\in\X_1$ for which the essential supremum is finite.

We define the following operators 
\begin{itemize}
\item The spatial \emph{mass} density of a velocity distribution,
\begin{equation}\label{mass}
\bar{}:\BV\ra \RR, \qquad \bar{p} = p(V). 
\end{equation}
Clearly, the operator \:$\bar{ }$\: is Lipschitz continuous.
\item The \emph{lifting} of a measure on $\Sn$ to a measure on $V$,
\begin{equation}\label{lifting} 
\tilde{}:{\mathcal B}(\Sn)\ra\BV, \qquad \tilde{q} = m \otimes q
\end{equation}

where $m$ is a probability measure on $[s_1,s_2]$. If $s_1=s_2$, then 
$\,\tilde{}\,$ just maps a measure on $\Sn$ to the same measure on 
$\{s_1\}\times\Sn$. In the paper \cite{Hillen06} it was taken to be the
normalized Lebesgue measure on $[s_1,s_2]$, which corresponds to the weight
parameter $\omega$ defined in \cite[equation (4)]{Hillen06}. The choice
$m([s_1,s_2])=1$ guarantees that 
\begin{equation*}
||\tilde{q}||_{L^\infty(\RRn,\BV)}  = ||q||_{\X_2}.
\end{equation*}
In particular, a function that takes values among the probability measures on
$\BS^+$ is mapped to a function taking values among probability measures on
$\BV^+$. Since $\,\tilde{}\,$ is a linear operator it is Lipschitz continuous.
Additionally, we use the lifting to connect the measures on $V$ and on $\Sn$ in
a natural way as
\begin{equation}\label{dvdtheta}
 \ud v = m\otimes \ud\theta,
 \end{equation}

\item The \emph{mean projection operator} (for undirected fibres) 
\begin{equation}\label{mean_projection}
\Pi_u(p)(\theta) =\frac{1}{\bar{p}} \int_V \left| \theta \cdot
\frac{v}{||v||}\right|\, \ud p(v).
\end{equation}
For sake of simpler notation and to avoid difficulties when $\bar{p}=0$, we
introduce the operator
\begin{equation*}
\Lambda:\X_1\ra L^1(\RRn,L^\infty(\Sn)), \qquad \Lambda(p) = \bar{p}\,\Pi_u(p).
\end{equation*}
Notice that $\Lambda$ is linear and if $||p||_\infty<\infty$ then
\begin{equation*}
||\Lambda(p)||_{L^\infty(\RRn,L^\infty(\Sn))}\le
||\bar{p}||_{L^\infty(\RRn,\RR)}.
\end{equation*}
For sake of completeness we also state the directed version of the operator
$\Lambda$,
\begin{equation*}
\Lambda_d(p)(\theta) = \int_V  \theta \cdot \frac{v}{||v||} \, \ud p(v).
\end{equation*}
As said above, existence of solutions is shown completely analogously in the two
cases.
\item The \emph{relative alignment operator} again, using the notation from
\cite{Hillen06}
\begin{equation}\label{relative_alignment}
A_u(p,q)=\int_{\Sn}\Pi_u(p)(\theta)\,\ud q(\theta).
\end{equation}
Similarly to the introduction of $\Lambda$, we will work with 
\begin{equation*}
B:\X \ra L^1(\RRn,\RR), \qquad B(p,q) = \bar{p}\,A_u(p,q).
\end{equation*}
Notice that $B$ is bilinear and if $||p||_\infty<\infty$, then
\begin{equation*}
||B(p,q)||_{L^\infty(\RRn,\RR)} \le
||\bar{p}||_{L^\infty(\RRn,\RR)}||q||_{\X_2}.
\end{equation*}
\end{itemize}
The operators $\Lambda$ and $B$ are Lipschitz continuous on bounded subsets. 

Let $\mu>0$ denote the turning rate and $\kappa>0$ denote the rate
of fibre degradation. The model \eqref{old_model} can be written as equality of
measures
\begin{equation}\label{Cauchy_problem}
\begin{aligned}
\frac{\partial p}{\partial t} +v\cdot \nabla p  &=-\mu p+\mu\bar{p}\,\tilde{q},
\\
\frac{\partial q}{\partial t} &= \kappa(\Lambda(p)-B(p,q)) q, \\
p(x,0)&=p_0(x), \quad q(x,0)=q_0(x).
\end{aligned}
\end{equation}

\section{Existence Results}\label{s:existence}

To provide a framework for local and global existence of solutions we define 
\begin{equation}\label{operator}
\begin{aligned}
D(\A)&=\{(p,q)\in\X\::\: \nabla p\in\X_1^n\, \}, \\
\A\left(\begin{array}{c}p\\q\end{array}\right) &=
\left(\begin{array}{cc} -v\cdot\nabla &
0\\0&0\end{array}\right)
\left(\begin{array}{c}p\\q\end{array}\right).
\end{aligned}
\end{equation}
Here $\nabla = \nabla_x$ is interpreted in the sense of weak derivatives of
Banach space-valued functions. We write $f=\nabla_x p$ for a function
$f\in\X_1^n$ if for all test functions $\phi\in W^{1,1}(\RRn,C(V))$
\begin{equation*}
-\int_{\RRn} f(x)\cdot\nabla_x \phi(x)\,\ud x =\int_{\RRn}p(x) \phi(x)\,\ud x
\in \BV,
\end{equation*}
where the integrals are Bochner integrals taking values in $\BV$.  Observe that
the domain $D(\A)$ is dense in $\X$, as it contains the space
$C^\infty(\RRn,\BV)\times \X_2$ of infinitely differentiable functions, which is
dense in $\X$ \cite[Theorem 2.16]{LiebLoss}. The operator $\A$ with domain
$D(\A)$ is the generator of a positive $C_0$-semigroup $\U(t)$ on the Banach
space $\X$ (see also Theorem 1 in \cite{Arlotti}).

Notice that the operator $-v\cdot\nabla$ is the \emph{collisionless transport
operator} occurring in the linear Boltzmann equation which has been studied by
many authors, see \cite{Greiner,Engel}, \cite[Chapter 13]{Kaper} and the
references therein. It generates a semigroup (in fact, a group) $\U_1$ on the
space $\X_1$ via
\begin{equation}\label{shifted_measure}
\U_1(t)p_0(x,A) =p_0(x-At,A):= \int_A  p_0(x-t\,\ud v,\ud v),
\end{equation}
for Borel sets $A\subset V$. Clearly, the positive cone $\X_1^+$ is invariant
under $\U_1$.  The group $\U_1$ preserves the $L^1$-norm while for
$||p_0||_\infty<\infty$ we have 
\begin{equation*}
|\U_1(t)p_0|(x,\cdot) = \int_V \ud p_0(x-t\,\ud v,\ud v) \le (1+ts_2|\Sn|)
||p_0||_\infty.
\end{equation*}
We denote the semigroup on $\X$ generated by the operator $\A$ from equation
\eqref{operator} by $(\U(t))_{t\ge0}$. It has a diagonal structure
\begin{equation}\label{def_semigrp}
\U = \left(\begin{array}{cc} \U_1 & 0\\0&I\end{array}\right),
\end{equation}
\begin{sloppypar}\noindent
where $I$ denotes the identity on $\X_2$. In the operator norm, $\U$ satisfies
\mbox{$||\U(t)||_{\mathcal{L}(\X)}\le 1$} and for $||u_0||_\infty <\infty$ we
obtain  
\end{sloppypar}
\begin{equation}\label{D_invar}
|| \U(t) u_0||_\infty \le (1+ts_2|\Sn|)  ||u_0||_\infty.
\end{equation}

For a pair $u=(p,q)\in\X$ define the map $||\,\cdot\,||_\infty:\X\ra[0,\infty]$
by
\begin{equation*}
||u||_\infty =  ||p||_\infty + ||q||_\infty,
\end{equation*}
and set
\begin{equation*}
\D = \{ u\in\X\::\:||u||_\infty  <\infty\}.
\end{equation*}
For every $r>0$ the set
\begin{equation*}
\D_r = \{ u\in\X\::\: ||u||_\infty \le r\}
\end{equation*}
is closed in $\X$ (in particular, the projection of $\D$ onto $\X_1$ is closed
with respect to the $L^1$-norm on $\X_1$). Indeed, let $p_n\in \X_1$ be a Cauchy
sequence with $||p_n||_\infty\le r$ for all $n$. Since $\X_1$ is complete, it
has a limit $p$. We claim that $||p||_\infty\le r$. Suppose that this were not
the case, then there would be an $\varepsilon>0$ and a set $A\subset\RR^n$ with
Lebesgue measure $|A|>0$ such that $||p(x)||_{\BV}\ge r+\varepsilon$ for all
$x\in A$. But then clearly the $L^1$-norm would satisfy
$||p_n-p||_{\X_1}\ge\varepsilon|A|>0$, which is a contradiction.

Problem \eqref{Cauchy_problem} can now be written as an abstract Cauchy problem
\begin{equation}\label{Cauchy_problem_abstract}
\begin{aligned}
u' &= \A u + F(u), \\
u(0)&=u_0,
\end{aligned}
\end{equation}
with $u=(p,q),\,u_0=(p_0,q_0)\in\D$. 
\begin{definition}\label{definition:mild_solution}
\cite{Oharu87} Let $u_0=(p_0,q_0)\in\D$. We say that a function $(p,q)=u\in
C([0,\infty),\D) $ is a \emph{global mild solution} if $F(u(\,\cdot\,))$ is
continuous and it satisfies the
integral equation
\begin{equation}\label{def_mild}
u(t) = \U(t)u_0+\int_0^t\U(t-s) F(u(s))\,\ud s,
\end{equation}
where $\U(t)$ is the  semigroup defined in equation \eqref{def_semigrp}. We call
a function $u=(p,q):[0,T)\ra\D$ a
\emph{classical solution} if it satisfies the following
properties
\begin{itemize}
    \item[(i)] $u\in C^1((0,T),\X)\cap C([0,T),D(\A))$, and
    \item[(ii)] equation \eqref{Cauchy_problem_abstract} holds.
\end{itemize}
\end{definition}
Our first result is
\begin{theorem}\label{theorem:1} Assume that $q_0(x,\Sn)=1$ for almost every
$x\in\RRn$, then the problem \eqref{Cauchy_problem_abstract} has a unique global
positive mild solution for every $u_0\in \D\cap\X^+$. 
\end{theorem}

\subsection{Proof of Theorem \ref{theorem:1}}\label{global}
The proof of Theorem \ref{theorem:1} is established in the following Lemmas.
\begin{sloppypar}
\begin{lemma}\label{lemma1} The right hand side of equation
\eqref{Cauchy_problem}
defines a nonlinear map \mbox{$F:\D\ra\D$}, which maps $\D$ into itself 
\begin{equation*}
F(p,q) = \left(\begin{array}{c}F_1(p,q)\\ F_2(p,q) \end{array}\right) =
 \left(\begin{array}{c} -\mu p +\mu\bar{p}\,\tilde{q} \\
\kappa(\Lambda(p)-B(p,q)) q \end{array}\right).
\end{equation*}
The map $F$ is Lipschitz continuous on bounded subsets of $\D$. 
\end{lemma}
\end{sloppypar}
\emph{Proof.} Observe that for $(p,q)\in\D$ the product $\bar{p}\,\tilde{q}$ is
well defined and  
\begin{equation*}
\begin{aligned}
||\bar{p}\,\tilde{q}||_{\X_1}&\le
||\bar{p}||_{L^1(\RRn,\RR)}||q||_{\X_2}=||p||_{\X_1}||q||_{\X_2}, \\
||\bar{p}\,\tilde{q}||_{L^\infty(\RRn,\BV)}&\le
||p||_{L^\infty(\RRn,\BV)}||q||_{\X_2},
\end{aligned}
\end{equation*}
in particular,
\begin{equation*}
||F_1(p,q)||_{L^\infty(\RRn,\BV)}\le 2\mu
||p||_{L^\infty(\RRn,\BV)}||q||_{\X_2}.
\end{equation*}
For functions $\varphi\in L^\infty(\Sn)$ and measures $q\in\BS$ we define the
product $\varphi q\in\BS$ by way of
\begin{equation}\label{function_times_measure}
(\varphi q)(M) =  \int_M\varphi(\theta)\,\ud q(\theta),
\end{equation}
where $M\subset \Sn$ is a Borel set. This multiplication extends to functions in
$L^\infty(\RRn\times\Sn)$ and $L^\infty(\RRn,\BS)$ and we have
\begin{equation*}
||\varphi q ||_{\X_2} \le  ||\varphi||_{L^\infty(\RRn\times\Sn)}||q||_{\X_2}.
\end{equation*}
With $\varphi(\theta)=\Lambda(p)(\theta)-B(p,q)$ we obtain
\begin{equation*}
||F_2(p,q) ||_\infty =||(\Lambda(p)-B(p,q)) q ||_{\X_2} \le
||\bar{p}||_{L^\infty(\RRn,\RR)}(1+||q||_{\X_2}) ||q||_{\X_2},
\end{equation*}
showing that $F_2$ takes values in $\X_2$. Computations similar to those just
carried out give the local Lipschitz continuity of $F$ on bounded subsets of
$\D$. For example, 
for $(p_1,q_1),\,(p_2,q_2)\in\D$ and
$||p_1||_{\X_1}+||p_2||_{\X_1}+||q_1||_{\X_2}+||q_2||_{\X_2} \le K$ there exists
 a constant $C(K)>0$ such that
\begin{equation*}
\begin{aligned}
||\bar{p}_1\tilde{q}_1-\bar{p}_2\tilde{q}_2||_{\X_1} &\le
||\bar{p}_1(\tilde{q}_1-\tilde{q}_2)||_{\X_1}+||(\bar{p}_1-\bar{p}_2)\tilde{q}
_2||_{\X_1}
\\
&\le C(||p_1-p_2||_{\X_1}+||q_1-q_2||_{\X_2}).
\end{aligned}
\end{equation*}
We omit the remaining calculations. \hfill $\Box$

\begin{lemma}\label{lemma2} Equation \eqref{Cauchy_problem_abstract} has a unique local
mild solution that remains positive for $u_0\in \X^+$.
\end{lemma}
\emph{Proof.} We set up a Banach's Fixed Point argument, but we cannot work on
$\D$  directly since that set is not complete. Hence we work with $D_R$ for some
$R$ large enough. For given $u_0\in \D$ and fixed $R,\,T>0$ we define
\begin{equation*}
E_{R,T}=\{u \in C([0,T],\D_R)\::\: u(0)=u_0\}.
\end{equation*} 
This set $E_{R,T}$ is a complete metric space, with the metric given by
\begin{equation*}
d(u,v) = \sup_{t\in[0,T]}||u(t)-v(t)||_\X.
\end{equation*} 
For  a function  $u\in E_{R,T}$  we define  
\begin{equation}\label{fixed_point_map}
\G u(t) = \U(t)u_0+\int_0^t\U(t-s) F(u(s))\,\ud s,
\end{equation}
this is again an element of $C([0,T],\D)$ with $\G u(0)=u_0$  (since $\D$ is
invariant under both the semigroup $\U$ and the nonlinearity $F$). 
We have for $u,\,v\in E_{R,T}$
\begin{equation*}
\begin{aligned}
||\G u(t)-\G v(t)||_{\X} &\le \int_0^t||\U(t-s)||_{\mathcal{L}(\X)}
||F(u(s))-F(v(s))||_\X \,\ud s \\
 &\le C t \sup_{s\in[0,t]}  ||u(s)-v(s)||_\X \le Ctd(u,v)
\end{aligned}
\end{equation*}
(where $C$ is the Lipschitz constant of $F$), hence 
\begin{equation*}
\begin{aligned}
d(\G u,\G v) \le CTd(u,v),
\end{aligned}
\end{equation*}
and by choosing $T$ sufficiently small, it can be achieved that $\G$ is a
contraction on the space $E_{R,T}$. If $v\in\D_R$, then we have (see the proof
of Lemma \ref{lemma1})
\begin{equation*}
\begin{aligned}
||F_1(v)||_\infty & \le 2\mu R^2, \quad ||F_2(v)||_\infty  \le R^2(1+R), \quad
\textrm{ and} \\
||F(v)||_\infty & \le  R^2(1+R+2\mu).
\end{aligned}
\end{equation*}
Let $u\in E_{R,T}$. We can estimate equation \eqref{fixed_point_map} 
\begin{equation*}
||\G u(t)||_\infty \le (1+ts_2|\Sn|)  ||u_0||_\infty +  t (1+ts_2|\Sn|) 
R^2(1+R+2\mu).
\end{equation*}
By choosing $R>2||u_0||_\infty$ (large) and $T$ small, namely  
\begin{equation*}
T \le \min\left\{\frac{1}{s_2|\Sn|},\, \frac{R-2 ||u_0||_\infty}{2
R^2(1+R+2\mu)}  \right\},  
\end{equation*}
we can achieve that 
\begin{equation*}
\sup_{t\in[0,T]}||\G u(t) ||_\infty \le ||\G u(T) ||_\infty< R,
\end{equation*}
for all $u\in E_{R,T}$. Hence the contraction $\G$ maps the complete metric
space $E_{R,T}$ into itself and hence has a unique fixed point by the Banach Fixed
Point Theorem. 

The positivity of solutions follows from the fact that the nonlinearity $F$ is
of multiplicative type. If either $p$ or $q$ becomes zero on a set at some time,
the left hand side of  equation  \eqref{Cauchy_problem_abstract} is
non-negative.
\hfill $\Box$

\medskip

Concerning the global existence of solutions, if $T_{max}(u_0)<\infty$, then by
\cite{Pazy}
\begin{equation*}
\lim_{t\nearrow T_{max}(u_0)} ||u(t)||_\X = \infty.
\end{equation*}
However, in our system \eqref{Cauchy_problem} a blow-up in finite time cannot
occur as the following lemma shows.
\begin{lemma}\label{lemma4}
Let $(p,q)(t)$ be a mild solution of equation \eqref{Cauchy_problem}
(equivalently, of  \eqref{Cauchy_problem_abstract}) taking
values in $\D\cap \X^+$. Then for all $t\in [0,T_{max})$ and almost every
$x\in\RRn$ we have
\begin{equation*}
q(x,t,\Sn)= 1, 
\end{equation*}
and there exists a constant $C>0$ such that  
\begin{equation*}
||p(t)||_{\X_1} =||p_0||_{\X_1}, \quad \textrm{and} \quad ||p(t)||_\infty \le
|V| ||p_0||_\infty e^{Ct}.
\end{equation*}
\end{lemma}
\emph{Proof.} Let $(p,q)$ be a mild solution of equation \eqref{Cauchy_problem}
taking values in $\X^+$. 
The second component of equation (\ref{def_mild}) reads 
\begin{equation*}
q(x,t) = I q_0(x) +\int_0^t I \kappa(\Lambda(p) - B(p,q) q)\, \ud s,
\end{equation*}
where $I$ denotes the identity. We evaluate this relation at $\Sn$ and use the
fact that 
\begin{equation*}
\begin{aligned}
\kappa(\Lambda(p)-B(p,q)) q(\Sn) &= \int_{\Sn} \Lambda(p)(\theta)\,\ud
q(x,\theta) -B(p,q)q(\Sn) \\&= B(p,q)(1-q(\Sn)).
\end{aligned}
\end{equation*}
We obtain
\begin{equation*}
1- q(x,t,\Sn) = 1-q_0(x,\Sn) - \kappa \int_0^t B(p,q)(1-q(x,s,\Sn)\, \ud s 
\end{equation*}
We apply Gronwall's lemma and obtain 
\begin{equation*}
1- q(x,t,\Sn)=(1-q_0(x,\Sn))\exp\left(-\kappa\int_0^tB(p(x,s),q(x,s))\,\ud
s \right).
\end{equation*}
The integrand is positive and bounded, hence by the assumption on $q_0$ we get 
\begin{equation}\label{probability}
q(x,t,\Sn)=1
\end{equation}
for almost all $x\in\RRn$. For the $L^1$-norm of $p$ we notice first that since
$p$ is positive, it satisfies
\begin{equation*}
p(x,t,V)=||p(x,t)||_{\BV}.
\end{equation*}
We evaluate the first equation of  \eqref{Cauchy_problem}  and obtain as a
consequence of \eqref{probability}
\begin{equation*}
\frac{\partial }{\partial t} \bar{p}(x,t) + \nabla\cdot\left( \int_Vv\,\ud
p(x,t,v)  \right) =0.
\end{equation*}
Integrating this equation over $\RRn$ gives
 \begin{equation*}
\frac{d }{dt}\int_{\RRn} \bar{p}(x,t)\,\ud x = -\int_{\RRn} \nabla\cdot\left(
\int_Vv\,\ud p(x,t,v)  \right)\,\ud x =0
\end{equation*}
by the divergence theorem. 
For the $L^\infty$-part of $p$ we use that $F_1(0,q) = 0$ and the following fact
\begin{equation*}
||F_1(p_1,q)-F_1(p_2,q)||_{L^\infty(\RRn,\BV)} \le \mu(1+||\tilde{q}||_{\X_2})
||p_1-p_2||_{L^\infty(\RRn,\BV)}.
\end{equation*}
We estimate from \eqref{def_mild} 
\begin{equation*}
\begin{aligned}
||p(t)||_\infty &\le |V|\left( ||p_0||_\infty + \int_0^t
||F_1(p(s),q(s))||_\infty\,\ud s\right)  \\ 
& = |V|\left( ||p_0||_\infty + \int_0^t
||F_1(p(s),q(s))-F_1(0,q(s))||_\infty\,\ud s \right) \\
&\le |V|\left(||p_0||_\infty + 
\mu(1+\sup_{s\in[0,t]}||\tilde{q}(s)||_{\X_2})\int_0^t||p(s)||_\infty\,\ud
s\right).
\end{aligned}
\end{equation*}
This inequality warrants application of Gronwall's lemma 
\begin{equation*}
||p(t)||_\infty \le |V| ||p_0||_\infty e^{Ct}
\end{equation*}
with a suitably chosen constant $C$.  By the density of the domain $D(\A)$ in
$\D\cap\X$ and because of the continuous dependence of the solution on the
initial datum we obtain the desired estimates for arbitrary initial data
$(p_0,q_0)\in\D\cap \X^+$.  \hfill $\Box$

Combining Lemmas \ref{lemma2}--\ref{lemma4} we conclude the proof of 
Theorem \ref{theorem:1}. 

\section{Steady States}\label{steady_states}
\begin{figure}[ht]
\begin{center}
\includegraphics[width=55mm]{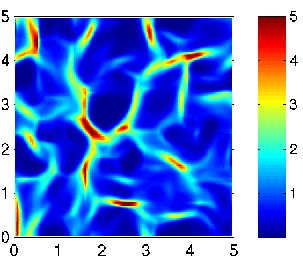}
\includegraphics[width=55mm]{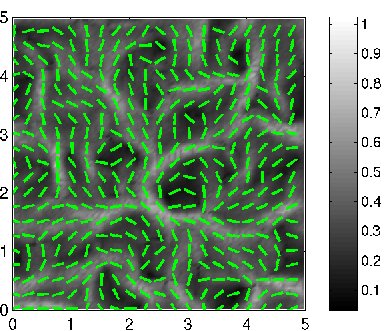}
\caption{Typical simulation of network formation for model
(\ref{Cauchy_problem}). The left figure shows the overall cell density $\bar
p(x,t)$ at a time where the steady state has almost been reached. Light (red)
color indicates high cell density and dark (blue) color indicates low cell
density. The figure on the right shows the underlying network, where the small
bars indicate the mean direction and the gray color describes the degree of
alignment. Light gray indicates highly aligned tissue, whereas dark gray/black
indicates close to  uniform distribution of directions. The simulations were
done by K.~Painter, and are described in detail in \cite{Painter}. We are
grateful to K.~Painter who allowed us to use this figure for illustrative
purposes.}\label{fig_Kevin}
\end{center}
\end{figure}

In numerical simulations by Painter \cite{Painter}, shown in  Figure
\ref{fig_Kevin}, we find interesting network patterns which form from random
initial data. Numerically, these patterns do not change after they have been
established. We expect that the system (\ref{Cauchy_problem}) allows for these
network patterns as steady states.  In this section we will develop a theory of
\emph{pointwise steady states} which are candidates for the observed network
patterns.

To describe steady states of (\ref{Cauchy_problem}) we introduce the bilinear
\emph{turning operator} 
\begin{equation*}
\LL\,:\BS\times\BV\ra\BV, \qquad \LL[q](p)= \tilde{q} \bar{p}-p.
\end{equation*}
Observe that in contrast to the paper of Chalub et al.~\cite{Chalub} the turning
kernel does not depend explicitly on  $v'$, i.e., the cells are reoriented
regardless of their original orientation. For $p\in\ker 
\LL[q]$, we have
\begin{equation*}
p = \tilde{q}\bar{p}.
\end{equation*}
Hence the orientation of the cells in a steady state is
entirely given by the fibre distribution $q$. This reflects the fact that a
perfect alignment of the cells with the underlying
 fibre network and only such a perfect alignment
remains invariant under the turning operator $\LL$.

The trivial steady state is a uniform distribution of fibres and cells:
\begin{lemma}\label{l:ss1}{\rm (Homogeneous tissue)} 
For every constant $\varrho\ge 0$ the pair 
\begin{equation*}
q(x) = \frac{\ud\theta}{|\Sn|},\qquad  p(x)= \varrho\tilde q = \varrho \frac{m
\otimes \ud\theta}{|\Sn|} = \varrho\frac{dv}{|\Sn|}
\end{equation*}
is a steady state of \eqref{Cauchy_problem} in $L^\infty(\RRn,\BV\times\BS)$.
The only steady state of this type in $\D\cap \X$ is obtained for $\varrho=0$.
\end{lemma}
\emph{Proof.}  If $q=\ud\theta/|\Sn|$ and $p= \varrho\tilde q$, then $\bar p =
\varrho$ and $p=\bar{p}\tilde{q}$. The right hand side of the first equation of
\eqref{Cauchy_problem} is zero. For the second equation, we need to compute
$\Lambda(p)-B(p,q)$. We have
\begin{equation*}
\Lambda(p)(x,\theta) =\varrho\int_V \left|\theta\cdot\frac{v}{||v||}\right|
\ud\tilde{q}(\theta) =\beta
\end{equation*}
for a  $\beta\ge0$, which is independent of $\theta$ and $x$. We obtain
\begin{equation*}
B(p,q)(x) =\int_{\Sn}\beta \frac{\ud\theta}{|\Sn|} = \beta.
\end{equation*}
Hence $\Lambda(p)-B(p,q)=0$ and the right hand side of 
the second equation of \eqref{Cauchy_problem} is zero as well.
Notice that the measures $\ud v $ and $\ud\theta$ are coupled in a natural way
through (\ref{dvdtheta}). 
\hfill$\Box$

\medskip
To find other steady states, we need a weak formulation.
\begin{definition}\label{definition:steady_state}
We say that $(p,q)\in \D\cap \X$ is a \emph{weak steady state} of
\eqref{Cauchy_problem}, if for each pair of test functions
\begin{equation*}
\phi \in W^{1,1}(\RRn, C(V)),\quad \psi\in C_0(\RRn,C(\Sn)), 
\end{equation*}
(where $C_0$ denotes functions vanishing at $\infty$) we have 
\begin{align}
&-\int_{\RRn}\int_Vv\cdot\nabla\phi(x,v)p(x,\ud v)\,\ud x=\nonumber \\
& \int_{\RRn}\int_V\phi(x,v)\left(-\mu p(x,\ud v)+\mu\bar{p}(x)\,\tilde q(x,\ud
v) \right)\,\ud x \label{ss1}, \\
&\int_{\RRn}\int_{\Sn} \left(\Lambda(p)(\theta)-B(p,q)\right)\psi(x,\theta)
q(x,\,\ud\theta)\, \ud x =0. \label{ss2}
\end{align}
\end{definition}
Notice that in this definition, $\phi$ and $\psi$ are real-valued functions in
the variables $v$ and $\theta$, respectively, hence the integrals on $V$ and
$\Sn$ make sense. In the next Lemma we study the biologically meaningful case of
a network completely aligned in a single direction. 

\begin{lemma}\label{l:ss2}{\rm (Strictly aligned tissue)}
Assume a preferred direction $\gamma\in\Sn$ is given and $\varrho\ge0$ is a
constant. Let $\delta_\gamma$ denote the Dirac mass on $\Sn$ concentrated at
$\gamma$. Then 
\begin{equation*}
p(x) =\varrho\tilde{q}, \qquad q(x) =\frac{\delta_\gamma+\delta_{-\gamma}}{2}
\end{equation*}
is a weak steady state  in $L^\infty(\RRn,\BV\times\BS)$. 
\end{lemma}
\emph{Proof.} Since $p\in\ker \LL$ and since it is spatially homogeneous,
equation \eqref{ss1} is satisfied. To study \eqref{ss2} we first compute the
following integrals on $\Sn$
\begin{equation*}
\begin{aligned}
\Lambda(p)(x,\theta) &= \frac{\varrho}{2} \int_V
\left|\theta\cdot\frac{v}{||v||}\right|
\,(\tilde{\delta}_\gamma+\tilde{\delta}_{-\gamma})(\ud v)=
\varrho|\theta\cdot\gamma|, \\
\int_{\Sn}\Lambda(p)(x,\theta) \psi(x,\theta)\,\ud q(\theta) &=
\frac{\varrho}{2}(\psi(x,\gamma)+\psi(x,-\gamma)) \\
B(p,q)(x) &= \frac{\varrho}{2}
\int_{\Sn}|\theta\cdot\gamma|\,(\delta_\gamma+\delta_{-\gamma})(\ud\theta) =
\varrho,
\end{aligned}
\end{equation*}
and
\begin{equation*}
\frac{B(p,q)(x)}{2}\int_{\Sn}\psi(x,\theta)\,(\delta_\gamma+\delta_{-\gamma}
)(\ud\theta) = \frac{\varrho}{2}(\psi(x,\gamma)+\psi(x,-\gamma)).
\end{equation*}
Hence 
\begin{equation*}
\int_{\Sn}\left(\Lambda(p)(x,\theta)-B(p,q)(x)\right)\psi(x,\theta)\,\ud
q(\theta) =0
\end{equation*}
for all $x\in\RRn$.
 \hfill$\Box$

\subsection{Pointwise Steady States}
In the preceding Lemmas we identified two simple homogeneous steady states. A
full analysis of other steady states at this level is difficult, since the very
weak formulation of measure valued solutions allows too many degrees of freedom.
We rather specialize to the study of \emph{pointwise steady states} as defined
below. With pointwise steady states, we can combine the previous two Lemmas and
design networks of aligned tissue with patches of uniform tissue. 

\begin{figure}
\begin{center}

\includegraphics[width=10cm]{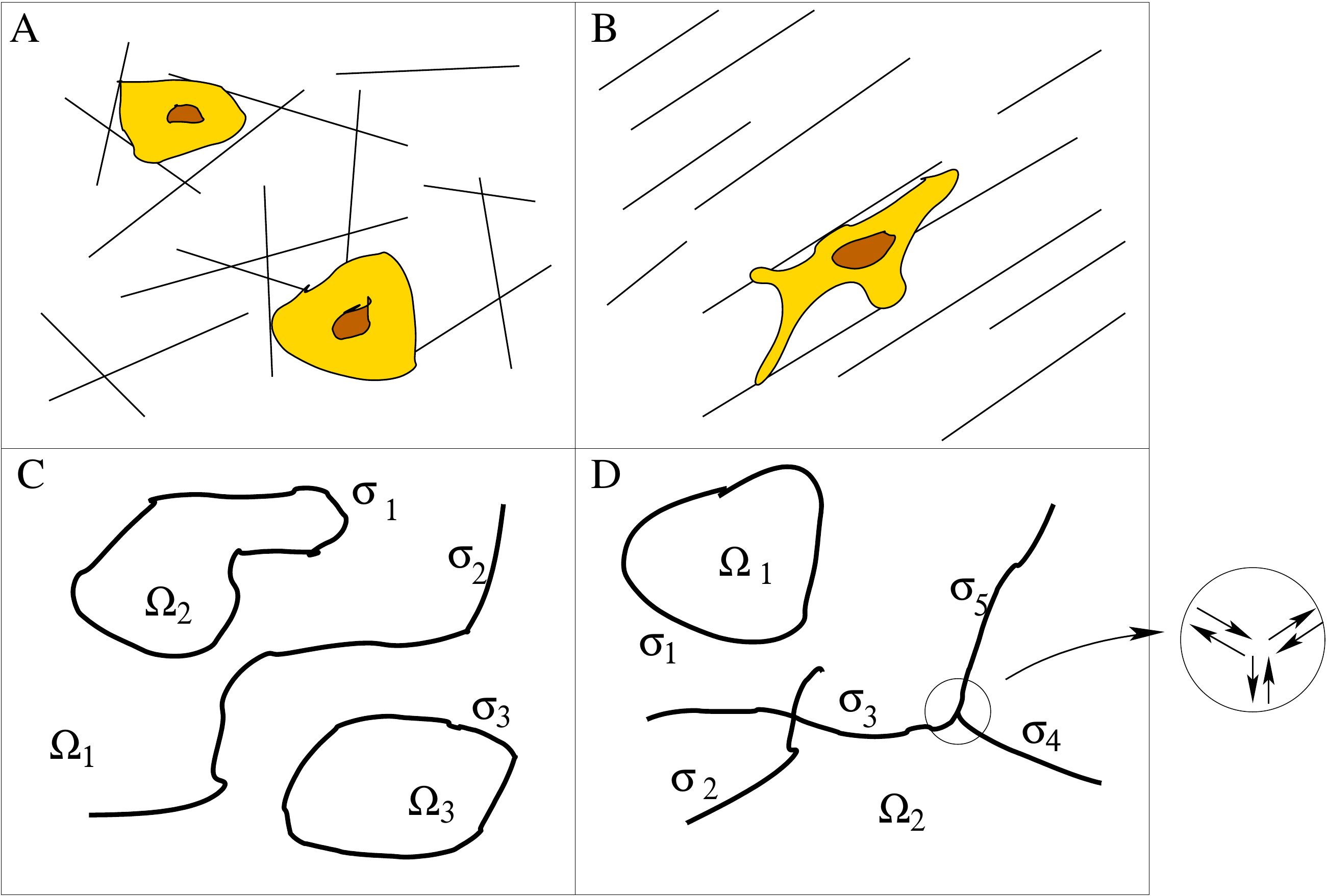}
\caption{Schematic of steady states. Figure A and B show typical fibre
alignment and cell alignment for the homogeneous steady state (in A) and the
strictly aligned steady state (in B). Figures C and D  show the geometric
construction that underlies
pointwise steady states. 
The fibres and cells are aligned tangentially along the curves $\sigma_i$ and
uniformly inside the domains $\Omega_i$. Figure C shows a pattern without
intersections, while Figure D shows intersections. One of the intersections has
been blown up to illustrate how a three pointed star of $120^\circ$ angles can
arise (see Corollary \ref{c:three}).}\label{f:steadystates}
\end{center}
\end{figure}

A schematic of the steady states which we construct here is 
shown in Figure \ref{f:steadystates}.


\begin{definition}\label{def-pointwise}{\rm
We say that $(p,q)\in \D\cap\X$ is a 
\emph{pointwise steady state} of \eqref{Cauchy_problem}, if  
\begin{enumerate}
\item $(p,q)$ is a weak steady state.
\item $p(x),\, q(x)$ is well defined for each $x\in\RRn$. 
\item For each test function $\Psi\in C(\Sn)$ and each $x\in\RRn$ 
\begin{equation}\label{def3}
\int_{\Sn} \left(\Lambda(p)(\theta)-B(p,q)\right)\Psi(\theta) q(x,\,\ud\theta)
=0.
\end{equation}
\item For each test function $\Phi\in C(V)$ and each $x\in\RRn$
\begin{equation}\label{def4}
\int_V\Phi(v)\left(-\mu p(x,\ud v)+\mu\bar{p}(x)\,\tilde q(x,\ud v) \right) =0. 
\end{equation}
\end{enumerate}
}  
\end{definition}

\begin{remark}{
\begin{description}
\item{(a)} An immediate consequence of item 1. and 4. in this definition is that
pointwise steady states satisfy
\begin{equation}\label{def5}
-\int_{\RRn}\int_Vv\cdot\nabla\phi(x,v)p(x,\ud v)\,\ud x= 0, 
\end{equation}
for each test function $\phi\in W^{1,1}(\RRn, C(V))$.
\item{(b)} Another immediate observation is that the homogeneous steady state
from Lemma \ref{l:ss1} and the completely aligned steady state from Lemma
\ref{l:ss2} are pointwise steady states. 
\end{description}
}\end{remark}

\medskip
In the following we classify pointwise steady states in $\RR^2$. It turns out
that the above definition allows for patchy steady states and for steady states
of network type. Patchy steady states include patches of uniform tissue
surrounded by areas of strictly aligned tissue, i.e.~a combination of the above
two types. Network type steady states arise if the areas of aligned tissue form
a connected network of curves with intersections and branches.
 For network type steady states, we will classify possible intersections of
network fibres. We are able to explicitly treat intersections of up to four
directions and we find a general algebraic condition for networks with 
intersections of higher order. 

\subsection{Patchy Steady States}
Assume a set of smooth curves $\sigma_i$, $i=1,\dots,N$ separate $\RR^2$ into
disjoint open sets $\Omega_i, i=1,\dots,k$.  Assume these curves $\sigma_i$ have
finite length, no intersections but they might be closed. Assume $p$ and $q$ are
uniform inside each patch 
\begin{equation}\label{pseins}
p(x)=p_i\frac{\ud v}{|V|}, \qquad q(x)=\frac{\ud\theta}{|\Se|}, \qquad \textrm{
if } x\in\Omega_i, \qquad i=1,\dots, k
\end{equation}
with  $p_i=0$ if $|\Omega_i|=\infty$. Since we are interested in pointwise
steady states, we need to define $(p(x), q(x))$
for $x\in\bigcup_i \sigma_i$. For each $i=1,\dots,N$ we denote the unit tangent
vector at $x\in\sigma_i$ by $\gamma_i(x)$, where we will suppress the argument
$x$ whenever possible. We define 
\begin{equation}\label{pszwei}
 q_i(x) = \frac{1}{2}(\delta_{-\gamma_i(x)}+\delta_{\gamma_i(x)}),\quad 
p(x) = \varrho_i \tilde q (x), \quad\mbox{for}\quad x\in\sigma_i.
\end{equation}
and $\varrho_i\geq 0$. 

\begin{lemma}\label{l:patchy}
The weak steady state defined by \eqref{pseins} and \eqref{pszwei} is a
pointwise steady state. 
\end{lemma} 
{\emph{Proof.} $(p(x),q(x))$ are defined for all $x\in\RR^2$ and 
as shown in the proofs of Lemma \ref{l:ss1} and \ref{l:ss2} the conditions 
(\ref{def3}) and (\ref{def4}) are satisfied for all $x\in \RR^2$. 
We only need to show that $(p,q)$ as defined above is a weak steady state, i.e.,
we need to confirm condition (\ref{def5}). 
 We find 
\begin{equation*}
\begin{aligned}
& -\int_{\RR^2}\int_Vv\cdot\nabla\phi(x,v)p(x,\ud v)\,\ud x\\ 
&= \sum_{i=1}^k \int_{\Omega_i}\int_V\phi(x,v)v\cdot\nabla \, p(x,\,\ud v) \,
\ud x - 
\sum_{i=1}^k\int_{\sigma_i}\int_V n\cdot v \,\phi(x,v)\, p(x,\,\ud v) \, \ud
\sigma,\\
&= 0 
\end{aligned}
\end{equation*}
The first integral vanishes since  $\nabla_x p(x,\,\ud v)=0$ 
in $\Omega_i$. The boundary integrals are zero, since we assumed that on
$\sigma_i$ the fibre orientation is tangential, i.e. $n\cdot v=0$ for all $v\in
\mbox{supp}\;p(x,t)$, where $n$ denotes the outer normal of $\sigma_i$ at
$x\in\sigma_i$.  
\hfill$\Box$

The above lemma allows for patches of uniform tissue surrounded by aligned
tissue. These could also be called encapsulations, as seen for many tumours in
tissue. A schematic of patchy steady states is given in Figure \ref{f:steadystates}C. 
The steady states in Lemma \ref{l:patchy} do, however, not allow for
intersections of the curves $\sigma_i$ so that they become of network type. To
obtain network steady states, we need to study possible intersections in more
detail. 
   
\subsection{Symmetric Intersections}
To study multiple directions we introduce two abbreviations. For a given vector
$\gamma\in \Se$ and a real valued function $\Psi$ on $\Se$ we define the
notation
\begin{equation*}
\delta_{|\gamma|} := \frac{1}{2}(\delta_{-\gamma}+\delta_{\gamma}), \quad
\Psi(|\gamma|) := \frac{1}{2}(\Psi(-\gamma)+\Psi(\gamma)).
\end{equation*}

We first consider the intersection of two directions $\gamma_1,\gamma_2\in \Se$
with $\gamma_1\neq\pm\gamma_2$ and with different weight $\alpha\in(0,1)$. For a
given $x\in\RR^2$ we define 
\begin{equation}\label{twointer}
 q(x) := \alpha\delta_{|\gamma_1|} +(1-\alpha)\delta_{|\gamma_2|}\quad
\mbox{and}\quad p(x) =\varrho\tilde q(x),
\end{equation}
where we set $\varrho=1$ without restriction. 
\begin{lemma}\label{l:twointer}
Assume $(p,q)$ is a weak steady state of (\ref{Cauchy_problem}) and at $x$ it is
of the form (\ref{twointer}). It can only be a pointwise steady state, if the
directions $\gamma_1$ and $\gamma_2$  have equal weight, i.e. if
$\alpha=\frac{1}{2}$.   
\end{lemma}
\emph{Proof.} We only need to check condition (\ref{def3}) of Definition
\ref{def-pointwise} at the intersection point $x$.
For this choice of $p$ and $q$ we find 
\begin{equation}\label{lambda1}
 \Lambda(p) =\int_V\left|\theta\cdot\frac{v}{||v||}\right|\ud\tilde
q=\alpha|\theta\cdot\gamma_1| +(1-\alpha)|\theta\cdot\gamma_2| 
\end{equation}
and
\begin{align}\label{b1}
B(p,q) &= \int_{\Se}\Lambda(p)(\theta) \ud q\nonumber\\
&=\alpha^2|\gamma_1\gamma_1| +2\alpha(1-\alpha)|\gamma_1\gamma_2|
+(1-\alpha)^2|\gamma_2\gamma_2|\nonumber\\
&= 2\alpha^2-2\alpha+1 +2\alpha(1-\alpha)|\gamma_1\gamma_2| 
\end{align}
where we used the fact that $|\gamma_i\gamma_i| =1, i=1,2$.  To check condition
(\ref{def3}), we need to test with a test function $\Psi\in C(\Se)$:
\begin{align*}
\int_{\Se}\Lambda(p)(\theta)\Psi(\theta)\, \ud q  &= \alpha (\alpha
+(1-\alpha)|\gamma_1\gamma_2|)\Psi(|\gamma_1|) \\
&\quad +(1-\alpha)(\alpha|\gamma_1\gamma_2| +1-\alpha)\Psi(|\gamma_2|),  \\
\int_{\Se}B(p,q)\Psi(\theta)\,\ud q &= B(p,q)\bigl(\alpha\Psi(|\gamma_1|)
+(1-\alpha)\Psi(|\gamma_2|)\bigr).\nonumber
\end{align*}
Hence to satisfy (\ref{def3}) for any test function, we need to satisfy 
\begin{equation}\label{cond2}
\alpha +(1-\alpha)|\gamma_1\gamma_2|= B(p,q) = \alpha|\gamma_1\gamma_2|
+1-\alpha.
\end{equation}
Comparing the first and last term, we obtain the equation
\[2\alpha-1 = (2\alpha-1)|\gamma_1\gamma_2|\]
which is satisfied only if $\alpha=\frac{1}{2}$. Notice that we assume
$|\gamma_1\gamma_2|\neq 1$. 
For $\alpha=\frac{1}{2}$ we find 
\begin{equation*}
B(p,q) = \frac{1}{2} +\frac{1}{2}|\gamma_1\gamma_2| 
\end{equation*}
and hence  the condition \eqref{cond2} is satisfied. \hfill$\Box$

Next we study the general case where at a given point $x\in\RR^2$ we have an
intersection of $N$-different directions $\gamma_1,\dots,\gamma_N\in\Se$. 
We study $N$ directions with equal weight: 
\begin{equation}\label{Ninter}
 q(x) := \frac{1}{N}(\delta_{|\gamma_1|} + \dots +\delta_{|\gamma_N|}) \quad
\mbox{and}\quad p(x) =\tilde q(x).
 \end{equation}
To decide if this intersection can be a pointwise steady state, 
we define a matrix of pairwise projections:
\begin{equation}\label{proj-matrix}
\Gamma :=(|\gamma_i \gamma_j|)_{i,j=1,\dots,N}.
\end{equation}
\begin{theorem}\label{t:G1}
Assume $(p,q)$ is a weak steady state of \eqref{Cauchy_problem}  and at $x$ it
is of the form \eqref{Ninter}. It can only be a pointwise steady state, if the
corresponding projection matrix $\Gamma$ has an eigenvector $(1,\dots,1)^T$.    
\end{theorem}
\emph{Proof.} Again, we only need to check condition \eqref{def3} of definition
\eqref{def-pointwise} at the intersection point. For the above choice of $p$ and
$q$ we find 
\begin{equation}\label{lambdaN}
 \Lambda(p) =\frac{1}{N}\sum_{i=1}^N|\theta \gamma_i| 
\end{equation}
and
\begin{equation}\label{b2}
B(p,q) = \frac{1}{N^2}\sum_{i,j=1}^N |\gamma_j\gamma_i|.
\end{equation}
Applied to a test function $\Psi\in C(\Se)$ we obtain
\begin{equation}\label{lambda2N}
\int_{\Se}\Lambda(p)(\theta)\Psi(\theta) \ud q  =
\frac{1}{N^2}\sum_{j=1}^N\left(\sum_{i=1}^N
|\gamma_j\gamma_i|\Psi(|\gamma_j|)\right) 
\end{equation}
and
\begin{equation}\label{b2N}
\int_{\Se}B(p,q)\Psi(\theta)\ud q = \frac{1}{N^2}\sum_{i,j=1}^N
|\gamma_j\gamma_i| \; \frac{1}{N}\sum_{k=1}^N\Psi(|\gamma_k|).
\end{equation}
To satisfy condition \eqref{def3} the right hand sides of \eqref{lambda2N} and
\eqref{b2N} have to coincide for any test function. In particular we need to
satisfy
\begin{equation}\label{condN}
\sum_{i=1}^N|\gamma_l\gamma_i| =\frac{1}{N}\sum_{i,j=1}^N
|\gamma_i\gamma_j|,\quad\mbox{for each} \quad l=1,\dots,N
\end{equation}
This condition implies
\begin{equation}\label{cond2N}
\sum_{i=1}^N|\gamma_l\gamma_i|=\sum_{i=1}^N|\gamma_k\gamma_i|\quad \mbox{for
each} \quad l,k=1,\dots,N.
\end{equation}
It can be directly verified that condition \eqref{cond2N} implies \eqref{condN}
and it also implies \eqref{lambda2N}=\eqref{b2N}. Hence \eqref{cond2N} is the
limiting condition. This condition implies that the row-sums of the matrix
$\Gamma$ are all identical, and since $\Gamma$ is a symmetric matrix, the column
sums also have the same value. In other words \eqref{cond2N} is equivalent with
the statement that $\Gamma$ has an eigenvector $(1,\dots,1)^T$. 
\hfill$\Box$

A schematic of a steady state with intersections is shown in Figure \ref{f:steadystates}D. 

\begin{remark} {\rm 
Notice that a related matrix to $\Gamma$ is  well known in linear algebra: the
Gram matrix 
\begin{equation*}
G = (\gamma_i \gamma_j)_{i,j}  
\end{equation*}
plays a role in coordinate transformations and the square root of the Gram
determinant is a measure for the volume element spanned by the vectors
$\gamma_1, \dots, \gamma_N$.  
}
\end{remark}

\medskip
\begin{example}\label{e:one}
As an example, we apply this general result to the two-directional case studied
in Lemma  \ref{l:twointer}. For two directions we have 
\begin{equation*}
\Gamma = \left(\begin{array}{cc} 1 & |\gamma_1\gamma_2|\\
|\gamma_1\gamma_2| & 1\end{array}\right),\quad\mbox{and}\quad
\Gamma\left(\begin{array}{c}1\\1\end{array}\right)=
(1+|\gamma_1\gamma_2|)\left(\begin{array}{c}1\\1\end{array}\right).
\end{equation*}
\end{example}
For three directions we obtain an interesting result:
\begin{corollary}\label{c:three}
Assume $(p,q)$ is a weak steady state of (\ref{Cauchy_problem}) and at $x$ it is
of the form (\ref{Ninter}) with $N=3$. It can only be a pointwise steady state,
if the three directions have equal angle, i.e. 
$|\gamma_1\gamma_2|=|\gamma_2\gamma_3|=|\gamma_3\gamma_1|$. 
\end{corollary}
\emph{Proof.} We use the criterion from Theorem \ref{t:G1}. The vector
$(1,1,1)^T$ is an eigenvector of $\Gamma$ if 
\begin{equation*}
1+|\gamma_1\gamma_2| + |\gamma_1\gamma_3| = 1 + |\gamma_1\gamma_2| +
|\gamma_2\gamma_3| = 1 + |\gamma_1\gamma_3| + |\gamma_2\gamma_3|
\end{equation*}
which implies
\begin{equation*}
|\gamma_1\gamma_2|=|\gamma_2\gamma_3|=|\gamma_3\gamma_1|.
\end{equation*}
\hfill$\Box$

A three pointed intersection has been illustrated in Figure \ref{f:steadystates}D.

\medskip
The classification of intersections of four directions is  a bit more complex.
\begin{corollary}\label{c:four}
Assume $(p,q)$ is a weak steady state of (\ref{Cauchy_problem}) and at $x$ it is
of the form (\ref{Ninter}) with $N=4$. It can only be a pointwise steady state,
if the pairwise equal angle condition (\ref{pairw}) is satisfied.  
\end{corollary}
\emph{Proof.}  To illustrate this case we introduce another abbreviation
\begin{equation*}
g_{ij} :=|\gamma_i\gamma_j|.
\end{equation*}
The corresponding projection matrix for four directions reads
\begin{equation*}
\Gamma = \left(\begin{array}{cccc}
1 & g_{12} & g_{13} & g_{14} \\
g_{12} & 1 & g_{23} & g_{24} \\
g_{13} & g_{23} & 1 & g_{34} \\
g_{14} & g_{24} & g_{34} & 1
\end{array}\right)
\end{equation*}
and the eigenvalue condition is given by 
\begin{align*}
1 + g_{12} + g_{13} + g_{14} &= 1 + g_{12} + g_{23} + g_{24} = 1 + g_{13} +
g_{23} + g_{34} \\ &= 1 + g_{14} + g_{24} + g_{34} .
\end{align*}
Hence we obtain six unknowns and three equations, which will not give us such a
complete solution as for three directions. We can, however, reduce the above
condition to a set of pairwise equal angle conditions
\begin{equation}\label{pairw}
|\gamma_1\gamma_2| = |\gamma_3\gamma_4|,\quad 
|\gamma_1\gamma_3| = |\gamma_2\gamma_4|,\quad
|\gamma_1\gamma_4| = |\gamma_2\gamma_3|.
\end{equation}   
If all angles are $\pi/2$ then this condition is satisfied.  \hfill$\Box$  

\subsection{Unsymmetrical  Intersections}
In the numerical simulations shown in Figure \ref{fig_Kevin} we observe
intersections that are asymmetric in the sense that for a direction $\gamma_1$
the opposite direction $-\gamma_1$ is not seen. Indeed, at an intersection point
$x$ various fibres come together. This means that at this point $x$ the network
has a number of directions $\gamma_1,\dots, \gamma_N$, but the opposite
directions are missing (it could of course happen by chance that
$\gamma_j=-\gamma_i$ for some $i,j$). Even though we assume that the
distribution function $q$ is symmetric almost everywhere, we find exceptional
points at those intersection points. In this section we show that unsymmetrical
intersections can arise as steady states in the framework developed here. 

Assume at a given point $x\in\RR^2$ we have an unsymmetrical intersection of
$N$-different directions $\gamma_1,\dots,\gamma_N\in\Se$ with equal weight: 
\begin{equation}\label{Ninternosym}
 q(x) := \frac{1}{N}(\delta_{\gamma_1} + \dots + \delta_{\gamma_N}) \quad
\mbox{and}\quad p(x) =\tilde q(x).
 \end{equation}
To decide if unsymmetrical intersections can arise as pointwise steady states,
we carry out the same computations as in the previous section. It turns out that
the computations change only marginally and we omit the details here. For
example formulas \eqref{lambda2N} and \eqref{b2N} use $\Psi(\gamma_i)$ instead
of $\Psi(|\gamma_i|)$. This implies that the conditions for their existence
remain the same. We summarize:
\begin{theorem}\label{t:G2}
\begin{enumerate}
\item
Assume $(p,q)$ is a weak steady state of \eqref{Cauchy_problem} and at $x$ it is
of the form \eqref{Ninternosym}. It can only be a pointwise steady state, if the
corresponding projection matrix $\Gamma$ has an eigenvector $(1,\dots,1)^T$.    
\item If $N=3$ then $(p,q)$ can only be a pointwise steady state, if the three
directions have equal angle, i.e. 
$|\gamma_1\gamma_2|=|\gamma_2\gamma_3|=|\gamma_3\gamma_1|$. 
\item If $N=4$ it can only be a pointwise steady state, if the pairwise equal
angle condition (\ref{pairw}) is satisfied.  
\end{enumerate}
\end{theorem}

\begin{remark}{\rm 
\begin{description}
\item{(a)} Indeed, unsymmetrical intersections of three directions with angles
of $120^\circ$ seem to be typical building blocks for the network shown in
Figure \ref{fig_Kevin}. 
\item{(b)} Although other intersections (symmetric or asymmetric) do exist
theoretically, they are rarely seen in simulations. This raises the question of
stability of these steady states. We defer this question to future studies. 
\end{description}}
\end{remark}

\subsection{Other Steady States}
We consider more general steady states where the cell distribution is a multiple
of the lifted fibre distribution, that is
\begin{equation*}
p(x) = \varrho(x)\tilde{q}(x),
\end{equation*}
where $\varrho \in L^\infty(\RRn)$ is the density of cells (or even $\varrho \in
L^1\cap L^\infty(\RRn)$). The minimal condition for such a pair to be a steady
state is
\begin{equation*}
\Lambda(p)(\theta) = B(p,q),
\end{equation*}
in particular, the left hand side is actually independent of $\theta$. Because
of the linearity of the operators $\Lambda$ and $B$, this condition becomes
\begin{equation*}
\varrho(x)\Lambda(\tilde{q}(x))(\theta) = \varrho(x) B(\tilde{q}(x),q(x)).
\end{equation*}
Wherever $\varrho \neq 0$ this condition can be stated as
\begin{equation}\label{condition_q}
\int_V \left| \theta \cdot \frac{v}{||v||}\right|\, \tilde{q}(x, \ud v) = 
\int_{\Sn}\int_V \left| \theta \cdot \frac{v}{||v||}\right|\, \tilde{q}(x, \ud
v) \,q(x,\ud \theta), 
\end{equation}
for almost all $x\in\RRn$. Now let us try the following ansatz
\begin{equation*}
q(x) = f(x) \delta_{|\gamma(x)|} + (1-f(x))\Sigma,
\end{equation*}
where $0\le f(x) \le 1$ and $\Sigma$ is the normalized Haar measure on $\Sn$.
Notice that even the predominant direction $\gamma$ may depend on $x$ at this
point. A calculation gives
\begin{equation*}
\begin{aligned}
\int_V \left| \theta \cdot \frac{v}{||v||}\right|\, \tilde{q}(x, \ud v)&=
\int_{\Sn}|\theta\cdot\psi|\, q(x, \ud \psi) \\
&= f(x) + (1-f(x)) \int_{\Sn}|\theta\cdot\psi|\, \Sigma(\ud\psi) =:C. 
\end{aligned}
\end{equation*}
The last term in the second line is independent of $\theta$ because of the
rotational invariance of $\Sigma$. Hence we have
\begin{equation*}
\int_{\Sn} C \,q(x,\ud \theta) = C \int_{\Sn} \,q(x,\ud \theta) = C.
\end{equation*}
Observe that 
\begin{equation*}
\overline{\varrho(x)\tilde{q}(x)} = \varrho(x)
\end{equation*}
and the right hand side of the first equation of \eqref{Cauchy_problem}
vanishes. Hence the condition for $\varrho$ is
\[
\varrho(x)\nabla\tilde{q}(x)+\tilde{q}(x)\nabla\varrho(x) = 0,
\]
for almost all $x\in\RRn$, 
which can be written as 
\begin{equation}\label{condition_rho} 
\nabla(\varrho(x)\tilde{q}(x)) = 0,
\end{equation}
In summary, if we have found $q$ that satisfies
\eqref{condition_q}, then $\varrho$ can be determined from the differential
equation \eqref{condition_rho}.

\section{Analysis For Given fibre Distribution}\label{characteristics}
 
In this section we assume that cells do not remodel the fibre network and that
$q(x,t)$ is a given distribution. The $p$-equation of system
\eqref{Cauchy_problem} has a simple structure for given $q$. For classical
solutions we can use the method of characteristics to find an explicit solution.
For given $v \in V$, the characteristic equation is $\frac{d}{dt}x(t)=v$. Hence
the  characteristic through $x_0\in\RRn$ is given by $x(t)=x_0+vt$. We can write
the first equation of \eqref{Cauchy_problem} as follows
\begin{equation}\label{char0}
\frac{d}{dt}p(x(t),t)+\mu p(x(t),t)=\mu \tilde{q}(x(t),t)\bar{p}(x(t),t),
\end{equation}

We evaluate equation \eqref{char0} at $V$ and obtain 
\begin{equation*}
\frac{d}{dt}\bar{p}(x(t),t)=\frac{d}{dt}p(x(t),t,V)=-\mu p(x(t),t,V)+\mu
\bar{p}(x(t),t)\tilde{q}(x(t),t,V)=0,
\end{equation*}
where we have used the fact that $\tilde{q}(x(t),t,V)=1$. Hence
$\bar{p}(x(t),t)$ is constant along characteristics. Equation \eqref{char0} is
equivalent to the equation 
\begin{equation*}
e^{-\mu t}\frac{d}{dt}(p(x(t),t)e^{\mu t})=\mu \tilde{q}(x(t),t)\bar{p}(x(t),t).
\end{equation*}
Integrating the above equation with respect to time, we obtain 
\begin{equation}\label{char1}
p(x(t),t)=e^{-\mu t}p(x_0,0)+\mu e^{-\mu t} \bar{p}(x(t),t) \int_0^t e^{\mu
s}\tilde{q}(x(s),s)\,\ud s.
\end{equation}
For a given $(x,t)\in \mathbb{R}^{n}\times \mathbb{R}^+$, we find the
anchor-point $x_0(v)=x-vt$ and the corresponding backward characteristic in the
direction $v$
\begin{equation*}
x(s)=x-vt+vs.
\end{equation*}
Applying this to equation \eqref{char1}, we have 
\begin{equation}\label{char2}
p(x,t)=e^{-\mu t}p_0(x-t \,\ud v)+ \mu \bar{p}(x,t)\int_0^t
e^{-\mu(t-s)}\tilde{q}(x-(t-s)\,\ud v,s)\,\ud s.
\end{equation}
This is an equality in the Banach space $\BV$ and the term $p_0(x-t\,\ud v)$ on
the right hand side has to be interpreted as the $v$-shifted measure  defined in
equation \eqref{shifted_measure}. The same notation applies to $\tilde{q}$. We
evaluate this measure $p(x,t)$ at $V$, i.e., we compute $\bar{p}(x,t)=p(x,t,V)$,
and obtain 
\begin{equation*}
\bar{p}(x,t)=e^{-\mu t}p_0(x-Vt,V)+ \mu \bar{p}(x,t)\int_0^t
e^{-\mu(t-s)}\tilde{q}(x-V(t-s),s,V)\,\ud s,
\end{equation*}
this is an equality between real numbers. The measure $\tilde{q}$ is
non-negative and for fixed $w\in V$ we have \mbox{$\tilde{q}(x-w(t-s),s)(V)=1$},
and
\begin{equation}\label{eqK}
K(x,t)=\mu \int_0^t e^{-\mu(t-s)}\tilde{q}(x-V(t-s),s,V)\,\ud s>0.
\end{equation}
Thus we get 
\begin{equation*}
(1-K(x,t))\bar{p}(x,t)=e^{-\mu t}p_0(x-Vt,V).
\end{equation*}
If $K(x,t)\ne 1$, then we can solve for $\bar{p}$ as 
\begin{equation}\label{pbar}
\bar{p}(x,t)=\frac{e^{-\mu t}}{1-K(x,t)}p_0(x-Vt,V).
\end{equation}
Then $\bar{p}$ can be used in \eqref{char2} to find an explicit solution
\begin{equation}\label{psol}
\begin{aligned}
p(x,t)=&e^{-\mu t}p_0(x-t\,\ud v)\\
&+\frac{\mu e^{-\mu t}}{1-K(x,t)}p_0(x-Vt,V) \int_0^t
e^{-\mu(t-s)}\tilde{q}(x-(t-s)\,\ud v,s)\,\ud s.
\end{aligned}
\end{equation}
Notice that this solution only depends on the initial condition $p_0$ and on the
fibre distribution $q$. To clarify, equation \eqref{psol} is again an equality
in $\BV$ and the right hand side is of the type ``measure +
number$\,\cdot\,$measure''.

Equation \eqref{eqK} simplifies  drastically in the special case of constant
fibre distribution $q(x,t)=q$. In this case, it follows that 
\begin{equation*}
K(x,t)=\mu \int_0^t e^{-\mu(t-s)}\tilde{q}(V)\,\ud s=\mu \int_0^t
e^{-\mu(t-s)}\,\ud s=1-e^{-\mu t}.
\end{equation*}
Equation \eqref{pbar} becomes 
\begin{equation}\label{star}
\bar{p}(x,t)=p_0(x-Vt,V),
\end{equation}

Equation \eqref{star} deserves some interpretation. $\bar{p}$ is the mass
density of particles of all velocities at point $(x,t)$, whereas $p_0(x-Vt,V)$
integrates the initial condition over the domain of dependence of the  point
$(x,t)$, the set $\{x-tv\::\:v\in V\}$. The velocity distribution at $(x,t)$
arises by following all characteristics through $(x,t)$ backwards (see Figure
\ref{fig1}). We call \eqref{star} a generalized Huygens principle. The solution
$p(x,t)$ from \eqref{psol} can then be written entirely in terms of the initial
condition
\begin{equation}\label{convex_combination}
p(x,t)=e^{-\mu t}p_0(x-t\,\ud v)+(1-e^{-\mu t})p_0(x-Vt,V)\tilde{q}.
\end{equation}
Using equation \eqref{star}, the explicit solution can also be written as
\begin{equation}\label{convexcombination2}
p(x,t)=e^{-\mu t}p_0(x-t\,\ud v)+(1-e^{-\mu t})\bar{p}(x,t)\tilde{q}.
\end{equation}
Hence the solution is a convex combination of the initial condition $p_0$ and
the current amount of cells $\bar{p}$ redistributed with respect to the
``controlling'' distribution $\tilde{q}$. We will use this observation in the
next section to rigorously prove convergence of the parabolic limit.

\begin{figure}[ht]
\begin{center}
\includegraphics[width=90mm]{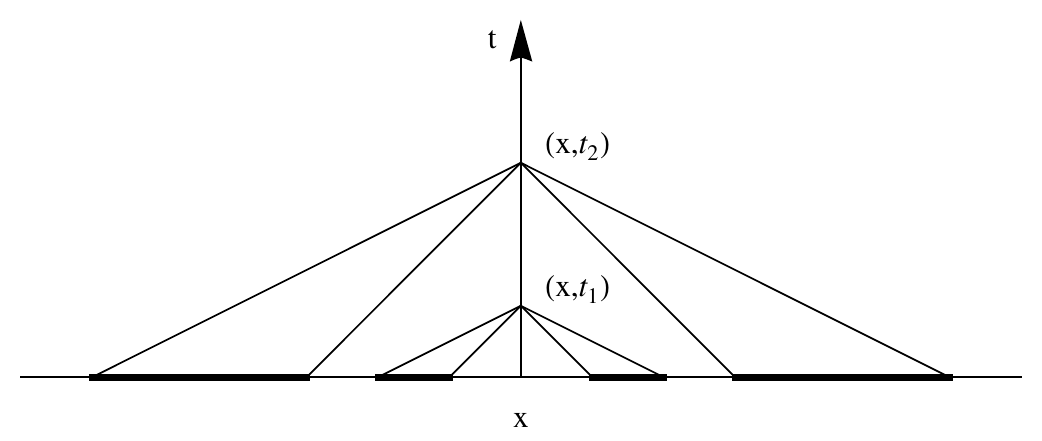}
\caption{The domain of dependence of the point $x$ at different time points is
shown as a thick solid line on the $x$-axis. In this example
$V=[s_1,s_2]\times\Sn $ with $s_1>0$ is an annulus.}\label{fig1}
\end{center}
\end{figure}

It is interesting to understand the biological meaning of these explicit
solutions. Equation \eqref{convexcombination2} tells us that eventually cells
will completely align to the given network structure. This has similarities with
glioma cell invasion along white matter tracks of the brain \cite{Jbabdi}.
Equation \eqref{psol} has a similar interpretation. In contrast to
\eqref{convexcombination2}, here the fibre distribution varies in time and
space. The integral terms in \eqref{psol} and \eqref{eqK} denote a temporal
average over the history of the fibre orientation, where the influence of the
history is exponentially damped. Then \eqref{psol} can be understood as cells
that try to align with the tissue while the tissue is changing. We see a
biological analogy with wound healing, where fibroblasts constantly modify the
collagen network while immune cells move through this fibre scaffold and heal
the wound \cite{Dallon:2001:MTE}. Notice that it is not our intention here to
model brain tumours or wound healing. These examples are only used as analogies
and thought experiments. It might be useful to make our model available to these
processes in the future.   
   
\subsection{The Parabolic Limit Problem}\label{limit}

As shown in \cite{Hillen06}, we can formally derive a diffusion limit equation
from equation \eqref{Cauchy_problem} under suitable scaling of space and time.
Let $\hat{x}$ and $\hat{v}$ denote reference length, and speed, respectively,
with the dimensionless quantity
\begin{equation*}
\e=\frac{\hat{v}}{\mu\hat{x}}
\end{equation*}
being small. We introduce rescaled variables as follows
\begin{equation*}
t^*=\e^2 t,\quad x^* = \frac{\e x}{\hat{v}}, \quad\textrm{and}\quad v^* = 
\frac{v}{\hat{v}}.
\end{equation*}
This gives
\begin{equation*}
\frac{\partial }{\partial t} = \e^2 \frac{\partial }{\partial t^*}, \quad 
\nabla_x = \frac{\e}{\hat{v}}\nabla_{x^*}, 
\end{equation*}
and we obtain, upon dropping the asterisks, the reduced parabolically scaled
equation 
\begin{equation}\label{eps_problem}\tag{P$_\e$}
\begin{aligned}
\e^2\frac{\partial p_\e}{\partial t} +\e v\cdot \nabla p_\e&=-\mu\LL[q](p_\e),
\\
p_\e(x,0)&=p_0(x) \in \D\cap\X_1^+.
\end{aligned}
\end{equation}
Simultaneously, we consider the limit problem
\begin{equation}\label{lim_problem}\tag{P$_0$}
\begin{aligned}
\frac{\partial \varrho}{\partial t} &=\nabla\cdot\left(D[q]\nabla\varrho\right),
\\
\varrho(x,0)&=p_0(x,V)=\bar{p}_0(x) \in L^{1,+}
(\RRn,\RR),
\end{aligned}
\end{equation}
with $||\varrho(\,\cdot\,,0)||_\infty<\infty$ and where the diffusion tensor is
given by
\begin{equation}\label{diffusion_tensor}
D[q] =  \frac{1}{\mu}\int_V v \otimes v \,\ud \tilde{q}(v).
\end{equation}
The formal derivation of the limit problem and the diffusion tensor in equation
\eqref{diffusion_tensor} was carried out in \cite[section 4]{Hillen06}, see in
particular equations (29) and (41) in that paper. We therefore omit these
calculations here. Notice that $D[q]$ can be written as the scaled
variance-covariance matrix $\mathbb{V}(q)$ of $q$, 
\begin{equation*} 
D[q] = \frac{1}{\mu}\int_{s_1}^{s_2}\int_{\Sn} (s\theta) \otimes (s\theta) \,\ud
q(\theta)\,\ud m(s)=\sigma \mathbb{V}(q),
\end{equation*}
where 
\begin{equation*} 
\sigma = \frac{1}{\mu}\int_{s_1}^{s_2} s^2\,\ud m(s), \qquad
\mathbb{V}(q)=\int_{\Sn} \theta\otimes\theta\,\ud q(\theta)
\end{equation*}
and we have used equation \eqref{lifting}.

We define the notion of a weak solution of equation \eqref{lim_problem}. 
\begin{definition}\label{definition:weak_solution}
Let $T>0$ be given. We say that $\varrho \in W^{1,1}([0,T],
\,W^{2,1}(\RRn,\RR))$ is a weak solution of \eqref{lim_problem} if the following
holds 
\begin{equation*}
\begin{aligned}
&-\int_{\RRn} \varrho(x,0)\phi(x,0)\,\ud x - \int_{\RRn}\int_0^T
\varrho(x,t)\frac{\partial \phi}{\partial t}(x,t) \,\ud t\, \ud x \\
&=-\int_{\RRn}\int_0^TD[q](x,t)\nabla \varrho(x,t)\cdot \nabla \phi(x,t) \,\ud
t\, \ud x
\end{aligned}
\end{equation*}
for all test functions $\phi\in C_c^\infty([0,T]\times\RRn)$ with
$\phi(\,\cdot\,,T)=0$, and in addition
\begin{equation*}
\varrho(x,0)=\bar{p}_0(x)
\end{equation*}
for almost every $x\in\RRn$. 
\end{definition}
The tensor $D[q]$ in equation \eqref{diffusion_tensor} is positive definite as
long as the support $supp\,q$ is not contained in a lower dimensional great
sphere. To see this we take $a\in\RRn$ and study
\begin{equation*}
a^T D[q] a = \frac{1}{\mu}\int_V (v \cdot a)^2 \,\ud \tilde{q}(v)> 0,
\end{equation*}
provided that $supp\,q$ is not contained in $\langle a\rangle^\perp \cap \Sn$
for any $a\in\RRn$. In this case, we have
the existence of weak solutions \cite{Friedman,Smoller}. 

\subsection{Convergence Result}\label{convergence_result}

The parabolic diffusion limit for chemotaxis was rigorously studied by Chalub
\textit{et al.} \cite{Chalub}. It was assumed that there exists a bounded
equilibrium velocity distribution $F(v)\in L^\infty(V)$ that is independent of
space, time and the distribution of the chemical signal, \cite[Assumption
(A0)]{Chalub}. The assumption (A0) in  \cite{Chalub} corresponds to our
assumption \eqref{Assumption_A0} below for the case that the equilibrium
distribution of the turning operator is a given function/measure on $\Sn$ (and
independent of $x$ and $t$). The difference arises from the fact that $F$ is
uniformly bounded while $\tilde q$ is a Borel measure.   

Since we are now equipped with a suitable functional analytical setting, we will
rigorously study the convergence to the parabolic limit. However, as shown
numerically in Painter (see \cite[Figure 9]{Painter}), the phenomenon of network
forming patterns is lost in the diffusion limit, hence we do not expect that
convergence to the diffusion limit is true in general. We assume that $q$ is
constant in space and time 
\begin{equation}\label{Assumption_A0}\tag{A0}
q(x,t) = q \in \BS
\end{equation}
for all $x\in\RRn$ and $t\ge0$ and that $q$ is symmetric with respect to $\theta
\mapsto -\theta$.
\begin{theorem}\label{theorem:2} 
Let assumption \eqref{Assumption_A0} hold and fix $T>0$. Let $(p_\e)_{\e\ge0}$
be the family of solutions to problem \eqref{eps_problem} and $\varrho$ the weak
solution to problem \eqref{lim_problem} (in the sense of Definition
\ref{definition:weak_solution}). Then, after possibly extracting a subsequence
we have the convergence
\begin{equation*}
p_\e \rightharpoonup \varrho \tilde{q}
\end{equation*}
in the weak${}^*$ topology on the space $L^\infty([0,T],\X_1)$. 
\end{theorem}
\emph{Proof.} Let $p_\e$ denote the $\BV$-valued solution of equation
\eqref{eps_problem}. We solve this equation as we did in Section
\ref{characteristics}, observing the new scaling with respect to $\e$. After
dividing equation \eqref{eps_problem} by $\e^2$ and applying
\eqref{convex_combination}, we find
\begin{equation}\label{alternative_21}
p_\e(x,t)=e^{-\frac{\mu}{\e^2} t}p_0\left(x-\frac{t\,\ud
v}{\e}\right)+(1-e^{-\frac{\mu}{\e^2}
t})p_0\left(x-\frac{Vt}{\e},V\right)\tilde{q}.
\end{equation}
The family $(p_\e)_{\e\ge 0}$ is uniformly bounded in $L^\infty([0,T],\X_1)$
since $||p(\,\cdot\,,t)||_{\X_1} = ||p_0||_{\X_1}$. 
Hence there exists a weak${}^*$-convergent subsequence, say $p_\e
\rightharpoonup p_*$ as $\e\ra0$.
Taking the $\BV$ norm in equation \eqref{alternative_21} and then taking the
supremum over all $x\in\RRn$ and $t\in [0,T]$ gives 
\begin{equation}\label{estimate_pe}
||p_\e(x,t)||_{\BV} \le \left\|p_0\left(x-\frac{t\,\ud
v}{\e}\right)\right\|_{\BV} +\left|p_0\left(x-\frac{Vt}{\e},V\right)\right| 
 \le 2||p_0||_\infty.
\end{equation}
and hence $||p_*(\,\cdot\,,t)||_\infty<\infty$. 
Using equation \eqref{star} in the rescaled coordinates we rewrite equation
\eqref{alternative_21} as
\begin{equation*}
p_\e(x,t)=e^{-\frac{\mu}{\e^2} t}\left( p_0\left(x-\frac{t\,\ud
v}{\e}\right)-\bar{p}_\e(x,t)\tilde{q}\right) + \bar{p}_\e(x,t)\tilde{q}.
\end{equation*}
Sending  $\e$ to $0$ in this equation we see that $p_*(x,t) = \varrho(x,t)
\tilde{q}$ for an appropriate function $\varrho\in L^1
(\RRn\times[0,T],\RR)$ with $||\varrho(\,\cdot\,,t)||_\infty<\infty$. It remains
to prove that $\varrho$ so defined satisfies the parabolic limit problem
\eqref{lim_problem}. To this end we define a residuum $r_\e$ and obtain with
\eqref{alternative_21}
\begin{equation*}
r_\e(x,t)=\frac{p_\e-\bar{p}_\e\tilde{q}}{\e} 
=\frac{e^{-\frac{\mu}{\e^2} t}}{\e}\left(p_0\left(x-\frac{t\,\ud
v}{\e}\right)-p_0 \left(x-\frac{Vt}{\e},V\right)\tilde{q}\right).
\end{equation*}
Observe that $\bar{r}_\e=0$ and for $\e\ge0$
\begin{equation*}
\frac{e^{-\frac{\mu}{\e^2} t}}{\e}\le 1.
\end{equation*}
By a similar argument as for $p_\e$ in \eqref{estimate_pe}, we get
\begin{equation*}
||r_\e(x,t)||_{\BV} \le  2||p_0||_\infty.
\end{equation*}
Hence there exists a  weak${}^*$-convergent subsequence $r_\e \rightharpoonup
r_*$. Finally, let $\varphi \in C_c^1(\RRn\times[0,T],\RR)$ be a test function
and observe that
\begin{equation*}
\begin{aligned}
&\e\int_0^T\int_{\RRn}\frac{\partial p_\e}{\partial t}(x,t)\varphi(x,t)\,\ud
x\,\ud t  \\
&=\e\int_{\RRn} p_\e(x,t)\varphi(x,t)\,\ud
x\Bigg|_0^T-\e\int_0^T\int_{\RRn}p_\e(x,t)\frac{\partial \varphi}{\partial
t}(x,t)\,\ud x\,\ud t.
\end{aligned}
\end{equation*}
Since the right hand side converges to zero as $\e\ra0$, so does the left hand
side and we have that
\begin{equation*}
\e\frac{\partial p_\e}{\partial t} \rightharpoonup 0
\end{equation*}
in the distributional sense. 

We divide equation \eqref{eps_problem} by $\e$ and obtain
\begin{equation*}
\e\frac{\partial p_\e}{\partial t} + v\cdot \nabla p_\e =-\mu r_\e.
\end{equation*}
Now we let $\e\ra0$, divide by $\mu$ and we obtain the following representation
of the limit of the residuum
\begin{equation}\label{r_infty}
r_* = -\frac{1}{\mu} v\cdot \nabla (\varrho \tilde{q}).
\end{equation}
We evaluate equation \eqref{eps_problem} at $V$ and obtain the conservation law
\begin{equation}\label{conservation_epsilon}
\e^2 \frac{\partial \bar{p}_\e}{\partial t} + \e \nabla\cdot\left(\int_V v \,\ud
(\e r_\e+\bar{p}_\e\tilde{q})\right) = 0.
\end{equation}
By the symmetry of $q$, we have
\begin{equation*}
\int_V v \,\ud \tilde{q}(v) = 0.
\end{equation*}
We divide equation \eqref{conservation_epsilon} by $\e^2$ and let $\e\ra 0$ and
obtain
\begin{equation*} 
\frac{\partial \varrho}{\partial t} =- \nabla\cdot \int_V v \,\ud r_*(v),
\end{equation*}
where
\begin{equation*} 
\frac{\partial \bar{p}_\e }{\partial t}  \rightharpoonup \frac{\partial
\varrho}{\partial t}
\end{equation*}
in the distributional sense. Using the representation \eqref{r_infty} we obtain
\begin{equation*}
\frac{\partial \varrho}{\partial t} =\nabla\cdot\left(\frac{1}{\mu}\int_V
v\otimes v \,\ud \tilde{q}(v)\,\nabla\varrho\right).
\end{equation*}
Hence $\varrho$ satisfies the limit equation \eqref{lim_problem}. \hfill$\Box$

\section{Discussion}\label{extension}
In this paper we consider mathematical properties of a model (\ref{old_model}), or equivalently 
\eqref{Cauchy_problem}, that describes mesenchymal cell movement in tissues. The
model was developed in \cite{Hillen06} and has been analyzed from various angles
in recent papers, \cite{Painter, Chauviere1, Chauviere2, WHL}. Through the
previous analysis it became evident that a solution framework is needed which
allows for measure valued solutions. Here we develop such a framework and prove
global existence of solutions in $\X$. We have used semigroup methods, since
they provide a dynamical systems point of view, and  we can use this framework
for linear stability analysis in future work. Alternative methods to show
existence include energy methods as developed by DiPerna and Lions
\cite{DiPerna}.

We were able to find non-trivial measure-valued steady states, which correspond
to homogeneous distributions, or to aligned tissue, or to patches of uniform
tissue with a network separating these patches. We found that pointwise steady
state show network properties as observed numerically. We also found that,
although the system has been formulated symmetrically, we can have unsymmetrical
intersection points. This confirms the 
interpretation in \cite{Hillen06}, where it was suggested that a network made
from undirected fibres can have characteristics of a directed network.  The
complete identification of steady states of \eqref{Cauchy_problem} is an
interesting open question. Furthermore it would be interesting to see whether
solutions of \eqref{Cauchy_problem} converge to steady states or to traveling
wave solutions as $t\ra\infty$. The existence result in $\X$ opens the door to a
rigorous linear stability analysis of steady states. This endeavor is left to
future work. 

The convergence to the parabolic limit is a standard feature of kinetic models
and it has been studied in many publications (see references in the text). Our
approach here extends known results to measure-valued solutions. Furthermore, we
formulate an explicit solution which shows that the solution basically is a
convex combination of initial data and its velocity-mean-value. The mean value
then is close to the parabolic limit. We also give an argument that the rigorous
convergence to a diffusion limit might only work for constant tissue.

Here we did not discuss the biological modelling of \eqref{Cauchy_problem}. We
would like, however, to discuss the biological
assumptions and propose various extensions, which could lead to more realistic
models.

One possibility is to introduce birth and death processes for the cells into the
model. For example, it is known that growth factors can be bound to the fibres
which promote proliferation. Also, harmful substances can be found in the fibre
network, possibly killing the cell. To model these effects, a term $G(p)$ would
be added to the right-hand side of the first equation of \eqref{Cauchy_problem}.
If $G(p)$ satisfies certain growth bounds, the global existence of solutions
will continue to hold.

A second possibility is to allow diffusion of $p$ with respect to both $x$ and
$v$. Cells are likely to undergo some random walk and may also change their
velocity randomly (a perfect alignment of cell velocities will disperse).
Diffusion with respect to the $x$ variable is easily modeled by adding a term of
the type $-D_x\Delta_x p$ to the left-hand side of the first equation. Also,
diffusion in the velocity can be modelled through an additional diffusion term
of the form $-D_v\Delta_v p$ (see also Dickinison \cite{Dickinson:2000:AGT} for
chemotaxis). For these cases we expect a smoothing property of the linear
semigroup and the totally aligned steady states will no longer exist.

Another possibility to expand and make the model more realistic would be to give
the fibres some elasticity and to let the fibres be moved by the cells. Also,
the cells should chose their new speed not randomly, but according to some
``stiffness'' of the neighborhood they are currently in. For example, a cell
that has to cut a lot of fibres in its way should slow down, while a cell that
is aligned well with the network can gain speed. Obviously, these are intuitive
ideas, and would have to be supported by biological evidence.

In model \eqref{Cauchy_problem} we implicitly assume that the protease is
released locally at the leading edge of the cell. In the literature, however,
various protease cutting mechanisms are discussed \cite{Friedl} and more detail
of the cutting can be included into the model (see also Painter \cite{Painter}). This might necessitate to
explicitly model the protease as a third variable through its own
reaction-diffusion equation.

A consideration of the length scales of the fibres relative to the size of the
moving cells might also give valuable input into the appropriate modelling
assumptions.

Finally, we have studied an unbounded domain $\RRn$ to avoid boundary
conditions. To formulate the correct boundary conditions for model
\eqref{Cauchy_problem} is not trivial. A commonly observed  effect seen in
tissue is that a tumour is encapsulated by a dense fibre network. In that case
the fibres at the boundary will be aligned tangentially to the boundary and
should trap moving cells inside the domain. The encapsulations can be understood
as patchy steady states, as described here.  A careful analysis of other
boundary conditions and its implications on existence and steady states is left
to future work. For the simulations in Figure \ref{fig_Kevin} K. Painter used
periodic boundary conditions on a square domain, i.e.~a flat torus.

\section*{\normalsize{Acknowledgments.}} We would like to thank Avner Friedman
(Ohio State University) for suggestions that led to the explicit solution in
Section \ref{characteristics} and Kevin Painter (Heriot-Watt University) for
providing Figure \ref{fig_Kevin}. We are greatly indebted to Pierre Magal
(University of Bordeaux) for enlightening discussions. The authors
thank the Institute for Mathematics and its Applications (IMA) at the University
of Minnesota, the Centre for Mathematical Biology (CMB) at the University of
Alberta, Edmonton and the CIRM Luminy, Marseille, where this research was done.

Received March 2009; revised March 2010.

\medskip


\begin{thebibliography}{99}

\bibitem{Alt80} (MR0661424)
   \newblock W. Alt,
   \newblock \emph{Biased random walk model for chemotaxis and related diffusion
approximation},
   \newblock J. Math. Biol. {\bf 9} (1980),  147--177.
\bibitem{Arlotti} (MR0899157)
   \newblock L. Arlotti,
   \newblock \emph{The {Cauchy} problem for the linear {Maxwell-Boltzmann}
equation},
   \newblock J. Diff. Eq. {\bf 69} (1987),   166--184.
\bibitem{Chalub} (MR2065025)
   \newblock F. A. C. C. Chalub, P. Markowich, B. Perthame, C. Schmeiser,
   \newblock \emph{Kinetic models for chemotaxis and their drift-diffusion
limits},
   \newblock Monatsh. Math. {\bf 142} (2000),  123--141.
\bibitem{Chauviere1} (MR2291824)
   \newblock A. Chauviere, T. Hillen, L. Preziosi,
   \newblock \emph{Modeling cell movement in anisotropic and heterogeneous
network tissues },
   \newblock Netw. Heterog. Media {\bf 2} (2007), 333--357.
\bibitem{Chauviere2} (MR2409219)
   \newblock A. Chauviere, T. Hillen, L. Preziosi,
   \newblock \emph{Modeling the motion of a cell population in the extracellular
matrix},
   \newblock Discr. Contin. Dyn. Sys. B {\bf Suppl. vol.} (2007), 250--259.
\bibitem{Cohn} (MR0578344)
   \newblock D. Cohn,
   \newblock ``Measure Theory'',
   \newblock Birkh{\"a}user, Boston, 1980.   
\bibitem{Dallon:2001:MTE}
   \newblock J. C. Dallon, J. A. Sherratt, P. K. Maini,
   \newblock \emph{Modelling the effects of transforming growth factor-$\beta$
on extracellular alignment in dermal wound repair},
   \newblock Wound Rep.\ Reg. {\bf 9} (2001), 278--286.  
\bibitem{Dickinson:2000:AGT} (MR1744041)
   \newblock  R. Dickinson,
   \newblock \emph{A generalized transport model for biased cell migration in an
anisotropic environment},
   \newblock J.\ Math. Biol. {\bf 40} (2000),  97--135.
\bibitem{DiPerna} (MR1014927)
   \newblock  R. J. DiPerna, P. L. Lions,
   \newblock \emph{On the {Cauchy} Problem for {Boltzmann} Equations: Global
Existence and Weak Stability },
   \newblock Ann. Math. {\bf 130} (1989),  321--366.
\bibitem{Engel} (MR1721989)
   \newblock K.-J. Engel, R. Nagel,
   \newblock ``One-Parameter Semigroups for Linear Evolution Equations'',
   \newblock Springer-Verlag, New York, Berlin, Heidelberg, 2000.
\bibitem{Friedl}
   \newblock P. Friedl, K. Wolf,
   \newblock \emph{Tumour cell invasion and migration: diversity and escape
mechanisms},
   \newblock Nat. Rev. Cancer {\bf 3} (2003), 362--374.
\bibitem{Friedman} (MR0181836)
   \newblock A. Friedman,
   \newblock ``Partial Differential Equations of Parabolic Type'',
   \newblock Prentice Hall, Englewood Cliffs NJ, 1964.
\bibitem{Greiner} (MR0731337)
   \newblock  G. Greiner,
   \newblock \emph{Spectral properties and asymptotic behavior of the linear
transport equation },
   \newblock Math. Z. {\bf 185} (1984),  751--775.
\bibitem{Hillen06} (MR2251791)
   \newblock T. Hillen,
   \newblock \emph{M$^5$ mesoscopic and macroscopic models for mesenchymal
motion},
   \newblock J. Math. Biol. {\bf 53} (2006), 585--616.
\bibitem{HKS1} (MR2139206)
   \newblock H. J. Hwang, K. Kang, A. Stevens,
   \newblock \emph{Global solutions of nonlinear transport equations for
chemosensitive movement},
   \newblock SIAM  J.  Math.  Anal. {\bf 36} (2005), 1177--1199.
\bibitem{HKS2} (MR2129381)
   \newblock H. J. Hwang, K. Kang, A. Stevens,
   \newblock \emph{Drift-diffusion limits of kinetic models for chemotaxis: a
generalization},
   \newblock Discrete Contin. Dyn. Sys. B. {\bf 5} (2005),  319--334.

\bibitem{Jbabdi} (MR000)
	\newblock A. Jbabdi, E. Mandonnet, H. Duffau, L. Capelle, K. R. Swanson,
M. Pelegrini-Issac, R. 
Guillevin, H. Benali, 
\newblock \emph{Simulation of Anisotropic Growth of Low-Grade Gliomas Using 
Diffusion Tensor Imaging}, 
\newblock Mang.\ Res.\ Med., {\bf 54} (2005), 616--624.
\bibitem{Kaper} (MR0685594)
   \newblock  H. G. Kaper,   C. G. Lekkerkerker, J. Hejtmanek,
   \newblock ``Spectral Methods in Linear Transport Theory'',
   \newblock Birkh{\"a}user, Basel, 1982.
\bibitem{LiebLoss} (MR1817225)
   \newblock  E. H. Lieb, M. Loss,
   \newblock ``Analysis'',
   \newblock American Mathematical Society, Providence RI, 2001.
\bibitem{Oharu87} (MR0883426)
   \newblock  S. Oharu,  T. Takahashi,
   \newblock \emph{Locally Lipschitz continuous perturbations of linear
dissipative operators and nonlinear semigroups },
   \newblock Proc. Amer. Math. Soc. {\bf  100} (1987),  187--194.
\bibitem{Painter} (MR2471301)
   \newblock  K. Painter,
   \newblock \emph{Modeling migration strategies in the extracelluar matrix},
   \newblock  J. Math. Biol. {\bf 58} (2009), 511--543. 
\bibitem{Pazy} (MR0710486)
   \newblock A. Pazy,
   \newblock ``Semigroups of Linear Operators and Applications to Partial
Differential Equations'',
   \newblock Springer-Verlag, New York, 1983.
\bibitem{Smoller} (MR1301779)
   \newblock J. Smoller,
   \newblock ``Shock Waves and Reaction-Diffusion Equations'',
   \newblock Springer-Verlag, New York, 1994.
\bibitem{Virga} (MR1369095)
   \newblock  E. G. Virga,
   \newblock ``Variational Theories for Liquid Crystals '',
   \newblock Chapman \& Hall, London, 1994.
\bibitem{WHL} (MR2452864)
   \newblock  Z. A. Wang, T. Hillen, M. Li
   \newblock \emph{Mesenchymal motion models in one dimension },
   \newblock SIAM J. Appl. Math. {\bf 69} (2008), 375--397.
\end{thebibliography}
\end{document}